\documentclass[12pt,twoside]{amsart}

\usepackage{geometry}
\usepackage{fancyhdr}

\usepackage{txfonts}
\usepackage{bbm}
\usepackage{mathrsfs}
\usepackage{amsfonts}
\usepackage{amsmath}
\usepackage{amssymb}
\usepackage{wasysym}
\usepackage{psfrag}
\usepackage{graphics,epsfig,amsmath}  
\usepackage{comment}
\usepackage{xcolor}
\usepackage{enumerate}
\usepackage{hyperref}
\usepackage{orcidlink}

\newtheorem{theorem}{Theorem}[section]
\renewcommand {\thetheorem}{\arabic{section}.\arabic{theorem}}

\newtheorem{proposition}{Proposition}
\renewcommand {\theproposition}{\arabic{section}.\arabic{proposition}}
\newtheorem{lemma}{Lemma}[section]
\renewcommand {\thelemma}{\arabic{section}.\arabic{lemma}}
\newtheorem{remark}{Remark}[section]
\renewcommand {\theremark}{\arabic{section}.\arabic{remark}}
\newtheorem{definition}{Definition}

\newcommand{\ZZ}{\mathbb{Z}}
\newcommand{\NN}{\mathbb{N}}
\newcommand{\RR}{\mathbb{R}}
\newcommand{\CC}{\mathbb{C}}
\newcommand{\TT}{\mathbb{T}}
\newcommand{\mc}{\mathcal}

\newcommand{\scr}{\mathscr}
\newcommand{\torus}{\mathbb{T}}
\newcommand{\real}{\mathbb{R}}

\allowdisplaybreaks

\numberwithin{equation}{section}

\title[foliation preserving torus map]{Resonances and Phase Locking phenomena for
foliation preserving torus maps}

\author[X. He]{Xiaolong He\orcidlink{0000-0002-3165-2883}}
\address{
School of Mathematics\\
Hangzhou Normal University\\
2318 Yuhangtang Rd.\\
Hangzhou, 311121, P.R. China}
\email{xlhe@hznu.edu.cn}
\thanks{Research of X.H.  supported by NSFC (NO. 12001148), ZJNSF (NO. LQ21A010013) and HZNU (NO. 2020QDL017)}

\author[R. de la Llave]{Rafael de la Llave\orcidlink{0000-0002-0286-6233}}
\address{
School of Mathematics\\
Georgia Institute of Technology \\
686 Cherry St. \\
Atlanta GA 30332, USA}
\email{rafael.delallave@math.gatech.edu}
\thanks{Research of R. L. supported in part by NSF DMS18000241}

\date{\today}

\begin{document}

\vspace{0.1in}

\begin{abstract}
It is well known for experts that resonances in nonlinear
systems lead to new invariant objects
that lead to new behaviors. 

The  goal of this paper is to  study the invariant sets  generated by resonances
under  foliation preserving  torus maps. That is
torus which preserve a foliation of irrational
lines $L_{\theta_{0}} = \{ \theta_0 + \Omega t |~ t \in \real \}
\subseteq \torus^{d}$. 

Foliation preserving maps  appear naturally
as reparametrization of linear flows in the torus and also play an
important role in several applications
involving coupled oscillators, delay equations, resonators with moving walls,
etc.
The invariant objects we find here, lead to predictions on the behavior
of these models.

Since the results of this paper are
meant to be applied for other problems, we have developed
very quantitative results giving very explicit descriptions of
the phenomena and the invariant objects that control them.

The
structure of the phase locking regions for  foliation preserving maps is
very different than for generic maps of the torus.
Indeed, for the sake of completeness, we have developed
similar analysis for the case of generic maps of the torus
and shown that the objects that appear in foliation
preserving maps are quantitatively and qualitatively different
from those of generic torus maps. This has consequences in applications.
\end{abstract}

\keywords{Foliation preserving torus map,  Phase locking,  Resonance,  KAM,  Lindstedt series,  Sternberg linearization,  Time reparameterization of irrational flows.
}

\subjclass[2010]{
  58F30, 
  37D10, 
  58F12. 
}
\maketitle

\section{Introduction}

The influential book \cite{Poincare} by Poincar\'{e} considered the study of
nearly integrable dynamical systems as \emph{the main problem of
mechanics}.

One importance of this study is that qualitatively new phenomena
may appear under some circumstances.
For example, when we consider perturbations of several oscillators, 
if the frequencies are rationally independent,
the perturbations average out and the system resembles for
a long time the  unperturbed system. On the other hand, if the
system presents a resonance (i.e. some of the frequencies
is a combination of others), the resonance can lead to
genuinely new phenomena not present in the original
system. For nonlinear systems, one needs to do several orders
of perturbation theory and new phenomena can
happen at each order of perturbation theory
and  there are phenomena  that happen beyond all orders of
perturbation theory.

In Hamiltonian systems, the geometry of the resonances
is extremely important for phenomena  such as Arnold diffusion
\cite{DH09, DdlLS06, Cor08} that show that perturbation theory
has limits.

In the  very non-resonant regions, one can continue averaging to all orders
and indeed the motion is similar to a rotation for all time. In
particular cases, one can apply KAM  theory and obtain
that indeed the system remains a rotation
(see  Section \ref{sec KAM}).


The goal of this paper is to study rather quantitatively, the
phenomena induced by resonances in some class of  torus maps with a
special structure (they preserve an irrational foliation).
See Subsection \ref{sect Foliation map} below for a precise
definition.  We also present KAM results that show that away from
resonances, indeed the system remains qualitatively the same.

The foliation preserving maps appear in several applications.
In pure mathematics, they appear as time reparametrizations of
constant vector fields in the torus. In more applied
contexts, they have appeared in the study of cavities with
moving boundaries (and, in general in quasi-periodic systems).
For us, a motivating application is the study of
state dependent delay equations
\cite{dlL21b-TBA}.

With a view to applications, we will present the results in
a very constructive way.
The new objects generated (or the KAM theory concluding that
the system remains a rotation) 
are obtained through solutions of a functional equation.
The perturbative analysis provides with approximate solutions
of the functional equation. The mathematical theory
includes also an a-posteriori result that shows that the
existence of approximate solutions implies the existence of
true solutions.  

Note that the usefulness of the a-posteriori
result goes beyond the perturbation regime  since the
approximate solutions of the functional equation can be produced by methods not based on
perturbation theory (e.g. by numerical calculations).
With a point of view to developing numerical methods, the
a-posteriori theorems are proved by establishing the
convergence of an iterative procedure. Implementing  the
iterative procedure provides an algorithm for the computation.

\subsection{Foliation preserving mappings}\label{sect Foliation map}
Let $\torus^d=\real^{d}/\mathbb{Z}^{d}$.
We say that $\Omega \in \real^d$ is an irrational frequency vector
when
$\Omega \cdot k  \ne 0$ for all $k \in \mathbb{Z}^d  \setminus \{0\}$.

If $\Omega$ is an irrational frequency vector, 
the sets
\begin{equation}
\label{leave}
L_{\theta_0} = \{ \theta_0 + \Omega t  \in \torus^d:~ t \in \real\}
\end{equation}
define a foliation of the torus $\torus^d$,
which we will denote by $\mathcal{F}_\Omega$.
Note that $L_{\theta_0}$ are the equivalent classes of
the equivalence relation
\[
\theta \sim \tilde \theta  \iff \exists~ t \in \real,~  s.t.~  \theta - \tilde \theta = \Omega t \ \text{mod}~  1.
\]
We are specially interested in maps $f$  of the torus which
preserve one of the above foliations (i.e., a leaf gets mapped into another leaf). More precisely, we are
interested in maps of the form:
\begin{equation} \label{foliationpreserving}
T_{\varphi}(\theta) = \theta + \varphi(\theta) \Omega,
\end{equation}
where $\varphi: \torus^d \rightarrow \real$
is a differentiable function.
To save notations, we write also $T_{\varphi}$ for the associated lift
 on $\real^{d}$.

With the notation \eqref{leave}, we
have that $T_{\varphi}( L_{\theta_0} )  \subseteq  L_{\theta_0}$.
Indeed, for any $\theta=\theta_{0}+\Omega t$ on $L_{\theta_{0}}$, there is
\begin{equation*}
  T_{\varphi} (\theta)=\theta_{0}+\Omega (t+\varphi (\theta_{0}+\Omega t))\in L_{\theta_{0}}.
\end{equation*}
Furthermore, let $\textrm{Diff}~(\torus^{d})$ be the set of  diffeomorphism mapping
on the torus and
let $\Xi$ be the subset of $\textrm{Diff}~(\torus^{d})$, in which the element $T_{\varphi}$ has
the form of $T_{\varphi}=Id+\varphi\Omega$.
Then one easily verifies that $\Xi$ is a subgroup of $\textrm{Diff}~(\torus^{d})$ under
the composition $T_{f}\circ T_{g}= T_{g+f\circ T_{g}}$
and the inverse of $T_{f}$ is given by
$T_{-f\circ T_{f}^{-1}}$.

Maps of the form \eqref{foliationpreserving}
appear naturally in the study  in resonant cavities affected
by quasi-periodic perturbations \cite{PdlLV03}, in  the study of
equilibria in quasi-periodic media \cite{dlLSZ16, dlLSZ17}
or in the study of state-dependent delay equations \cite{HdlL16, HdlL16b}.
In these applications, the objects we describe have
direct consequences: in the problems of  cavities they lead
to exponential growth of energy \cite{PdlLV03} and in state dependent delay equations,
they lead to lack of analyticity.
On the purely  mathematical side,
they appear as reparametrization of linear flows of the torus
\cite{ FKW01, Fay02, Cor02}.

In contrast with the generic torus maps, the special structure
of \eqref{foliationpreserving} leads to rather different
dynamical consequences.
For instance, there exists at most one non-vanishing Lyapunov
exponent for the foliation preserving torus maps.
More precisely, we have the following result, whose proof
can be found in \cite[Proposition 3.1]{PdlLV03}.

\begin{proposition}\label{LyaExp}
  Let $T_{\varphi}$ be a $C^{1}$ torus map
  of the form \eqref{foliationpreserving}.
  Then for every $\theta\in\torus^{d}$, $d-1$ of the Lyapunov exponents of
  $T_{\varphi}$ are zero. Besides these $d-1$ trivial exponents, for
  almost every $\theta$ $($in the sense of  any $T_{\varphi}$-invariant measure$)$, there is
  one Lyapunov exponent corresponding to the direction of $\Omega$.
\end{proposition}

Furthermore, if the
function $\varphi$ in \eqref{foliationpreserving} is bounded away from zero,
then there is no periodic points for $T_{\varphi}$. Indeed, since
the line $L_{\theta}=\{\theta+\Omega t: t\in\real\}$ in the cover
is mapped to itself,  the motion
on the line is increasing, so that no orbit can come
back to itself. So, there are
no periodic points in the line.
From the irrationality of $\Omega$, we obtain that
two different points on the line are also two different points in $\torus^{d}$,
so that all points in the orbit are different.

Periodic points in maps  of the form \eqref{foliationpreserving}
can exist but only when $\varphi$ has different signs. Clearly
if $\varphi(p) = 0$, the point $p$ is a fixed point.

\subsection{Resonant frequency}\label{resonant intrinsic}

We say that  $\Omega\in\RR^{d}$ is \emph{resonant} if there exist
$k\in\ZZ^{d}\setminus\{0\}$ and
$n\in\ZZ$ such that $\Omega\cdot k-n=0$. The frequency $\Omega\in\RR^{d}$ is
\emph{non-resonant} if for any $k\in\ZZ^{d}$ and $n\in\ZZ$, the relationship
$k\cdot\Omega-n=0$ implies $k=0$ and $n=0$.

Note that non-resonant, implies irrational because $0$ is a
particular case of an integer. However, it could happen
that for an irrational $\Omega$ we have 
that for some $k \in \mathbb{Z}^d \setminus\{0\}$ there is
 $k\cdot \Omega  = n $, where $n$ is  non-zero integer.
In such a case, the $\Omega $ could be irrational but resonant.

Denoting
\begin{equation}\label{resonant module}
\Gamma_{\Omega}(\ZZ)=\{k\in\ZZ^{d}|\  k\cdot \Omega \in \ZZ\},
\end{equation}
we see that if $k,\tilde{k}\in\Gamma_{\Omega}(\ZZ)$ so is $k+\tilde{k}$.
In more algebraic  language,  we note that 
$
\Gamma_{\Omega}(\ZZ)$ is a $\ZZ$-module and called the \emph{resonance module}
of $\Omega$.

Assume $\Omega$ is resonant (i.e. $\Gamma_{\Omega}(\ZZ)\neq\{0\}$), then there exist $k_{1},\cdots,k_{d-r}\in\ZZ^{d}
\setminus\{0\}$ linearly independent over $\RR$  such that $\{k_{j}\}_{j=1}^{d-r}$ form a
basis of $\Gamma_{\Omega}(\ZZ)$, which means
\begin{equation*}
 \Gamma_{\Omega}(\ZZ)=\left\{z\in\ZZ^{d}:\ z=\sum_{j=1}^{d-r}t_{j} k_{j}\ \textrm{with}\ t_{j}\in\ZZ \
\textrm{uniquely determined}\right\}.
\end{equation*}
Moreover, we can also find  a matrix $\mathfrak{A}\in\mathrm{SL}(d,\mathbb{Z})$,
$\omega\in\real^{r}$ and $L\in\mathbb{Z}^{d}$ in such a way that
\begin{equation}\label{intrinsic fre}
 \mathfrak{A}\Omega=
\begin{pmatrix}
 \omega\\
0
\end{pmatrix}
+L
\end{equation}
with $\omega\cdot m\neq 0$ for $m\in\ZZ^{r}\setminus
\{0\}$.

We refer to $\omega$ in \eqref{intrinsic fre} as the \emph{intrinsic frequency} of
$\Omega$. It is essentially unique, i.e., unique up to
change of basis in $\RR^{r}$ given by a matrix in $\mathrm{SL}(r,\ZZ)$.

\medskip
In our case, it is important to realize
that, given an  irrational frequency $\Omega'$, it is possible
that there  exists a real number $\alpha$ such that
$\alpha\Omega'$ is resonant. However, in this case, the resonance module of $\alpha\Omega'$ is only
one dimension (see also \cite{dlLSZ16}).

\begin{proposition}\label{one dimension resonance}
Let $\Omega'\in\real^{d}$ be irrational  and $0\neq\alpha\in\real$. If $\alpha\Omega'$ is resonant, then
the dimension of $\Gamma_{\alpha\Omega'}(\ZZ)$ is exactly one.
\end{proposition}

\noindent\textbf{Proof.}
For $k_{1}\cdot\alpha\Omega'-n_{1}=0$ and
$k_{2}\cdot \alpha\Omega'-n_{2}=0$, we see that
\begin{equation*}
 \alpha=\frac{n_{1}}{k_{1}\cdot\Omega'}=\frac{n_{2}}{k_{2}\cdot\Omega'}.
\end{equation*}
The irrationality of $\Omega'$ implies that $n_{2}k_{1}-n_{1}k_{2}=0$, which proves the proposition.
\qed

Note that for frequencies $\Omega$ that have a one-dimensional resonance
module, the intrinsic frequency has one dimension less.

Note also that if $\Omega'$ is irrational and $\alpha$ real, we have
that $\alpha \Omega'$ is irrational. However if we pick
$k \in \mathbb{Z}^d \setminus 0$ and $n \in \mathbb{Z}$,
if we take $\alpha = n/( k \cdot \Omega')$, we
have $k\cdot \alpha \Omega' = n$. Therefore the set  of $\alpha$ 
such that  $\alpha \Omega' $  is resonant is dense in
$\mathbb{R}$.

An important type of non-resonant frequency that will be used in this manuscript is
the following.

\begin{definition}\label{Diophantine}
We say that $\omega \in\RR^{s}$ satisfies the Diophantine condition of type $(\nu,\tau)$ if
 \begin{equation}\label{Diophantine condition}
  | k\cdot \omega -n|\geq \nu|k|^{-\tau}
 \end{equation}
for all $k\in\ZZ^{s}\setminus\{0\}$ and
$n\in\ZZ$ .
\end{definition}

We denote $\mathscr{D}_{s}(\nu,\tau)$ the set of vectors in $\RR^{s}$ satisfying the Diophantine conditions of type
$(\nu,\tau)$. It is well known that the set $\mathscr{D}_{s}(\tau)=\cup_{\nu>0} \mathscr{D}_{s}(\nu,\tau)$
occupies full Lebesgue measure in any open sets of $\RR^{s}$ for $\tau>s$. See \cite{dlL01}.

\subsection{Preliminaries}

In this paper, we deal with analytic functions and finitely differentiable functions defined
on the torus. We collect here some standard notations and results which
appear frequently later.

For the finite dimension space $\RR^{s}$ we adopt its supremum norm.
We  denote by $\mathscr{A}_{\rho}=\mathscr{A}_{\rho}(\TT^{d},\RR^{s})$
 the set of real analytic periodic functions which
are analytic in the complex neighborhood $\TT_{\rho}^{d}$ of $\TT^{d}$ in the complex space.
Here $\TT_{\rho}^{d}=\{x\in\CC^{d}/\ZZ^{d}: |\textrm{Im}x|\leq \rho\}$. We also endow
$\mathscr{A}_{\rho}(\TT^{d},\RR^{s})$  the supremum norm defined by
\begin{equation}\label{analytic norm}
\|w\|_{\rho}=\sup\{|w(x)|:\ x\in \TT_{\rho}^{d}\}.
\end{equation}

Then by the Cauchy estimate, it is readily seen that, if $w\in\mathscr{A}_{\rho}$,
the partial derivative with respect to its $j$-th argument $x_{j}$ satisfies
$$
\|w_{x_{j}}\|_{\rho-\sigma}\leq \sigma^{-1} \|w\|_{\rho}
$$
for all $0<\sigma<\rho$ and $1\leq j\leq d$. Furthermore, $w\in\mathscr{A}_{\rho}$
can be expanded into Fourier series
$$
w(x)=\sum_{k\in\ZZ^{d}}\widehat{w}(k) e^{2\pi i \langle k, x\rangle},
$$
whose Fourier coefficients $\widehat{w}(k)$ satisfies
\begin{equation}\label{exponential decay}
| \widehat{w}(k)|\leq \|w\|_{\rho}e^{-2\pi |k|\rho}.
\end{equation}

Next, we denote by $\mathscr{C}^{n}(\TT^{d},\RR^{s})$ the Banach space of $n$ times continuously
differentiable periodic functions, whose norm is given by
\begin{equation}\label{norm c n}
\|w\|_{\mathscr{C}^{0}}=\sup_{x\in\TT^{d}}|w(x)|,\quad
\|w\|_{\mathscr{C}^{n}}=\max_{0\leq j\leq n}\|D^{j}w\|_{\mathscr{C}^{0}}.
\end{equation}
If $w\in\scr{C}^{n}(\TT^{d},\RR^{s})$, we have that
\begin{equation}\label{poly decay}
 \sup_{k\in\ZZ^{d}}(|\widehat{w}(k)|\cdot |k|^n)\leq C_{n}\|w\|_{\scr{C}^{n}},
\end{equation}
where $C_n$ is the constant that depends only on $n$.
For more details on the regularities of functions characterized by the decay rate of Fourier
coefficients, we refer to \cite{dlL01, Gra14}.

Finally, we collect the following lemma on the  cohomology equation
\eqref{cohomology-torus-1} frequently used in KAM theory, whose proof can be
found in \cite{Rus75, dlL01}. We will encounter
 equation \eqref{cohomology-torus-1} when developing the resonant normal form,
formulating Lindstedt series near resonance and proving the KAM theorem for the foliation preserving
torus map.

\begin{lemma}\label{cohomology-torus}
Assume that $\omega$ satisfies the Diophantine condition \eqref{Diophantine condition}.
Let $Q\in\scr{A}_{\rho}$  be periodic function
with zero average, i.e., $\int_{\TT^{d}} Q~ \mathrm{d} \theta=0$.
Then there is a unique  solution $W$ of
\begin{equation}\label{cohomology-torus-1}
W(\theta)-W(\theta+\omega)=Q(\theta)
\end{equation}
such that $W$ has zero average.
Moreover, we have for all $0<\sigma<\rho$,
\begin{equation}\label{cohomo-esti}
\|W\|_{\rho-\sigma}\leq C
\nu^{-1} \sigma^{-\tau}\|Q\|_{\rho},
\end{equation}
where the constant $C$ depends only on the
Diophantine exponent $\tau$ and the dimension of the space.
\end{lemma}

\subsection{Organization of this paper}
The left of the paper is organized as follows. 
In Section \ref{sec r}, we study the
resonance and phase locking phenomena for the foliation preserving torus maps (FPTM).
We first establish a resonant normal form for FPTM, based on which an invariant surface
is constructed under some non-degeneracy conditions (see Theorem \ref{phase locking theorem}). Next we compute the Lindstedt series
for the invariant surface, and show the relationship between the Lindstedt series and 
the true embedding of the invariant surface (see Theorem \ref{Lind theorem}).
For the  two dimensional  FPTM, we study more closely the dynamics around the invariant
circle. For instance, we establish an extension of the Sternberg linearization theorem for the invariant circle (See
Theorem \ref{Sternberg} and \ref{constant expanding}), which
reduces the  dynamics of the invariant circle in the normal direction to the linear flow.
Moreover, we  study the structural stability of FPTM in two dimension.
To illustrate our theoretical results, we present a simple example at the end of Section \ref{sec r}.

In Section \ref{sec KAM}, we establish an a-posterior  KAM theorem for the
foliation preserving maps, which is devoted to the conjugation  problem of FPTM to the rigid rotation.
  In  Appendix \ref{general torus map}, we study the phase locking phenomena
for the general torus maps in the perturbative setting, which provides a comparison
to those special maps preserving the foliation and might be of interest itself.

\section{Resonance and phase locking phenomena}\label{sec r}

In this section, we study the phase locking phenomena for a family of foliation
preserving torus maps
close to  rotations on the torus.
 More precisely, we study
the following families of torus maps
\begin{equation}\label{original system}
 F_{\alpha, \varepsilon}(x)=x+\alpha\Omega+\varepsilon f_{\varepsilon}(x)\Omega,
\end{equation}
where $x\in\torus^{d}$, $\alpha\in\RR$, $\Omega\in\RR^{d}$ is irrational
and $f_{\varepsilon}$ is  defined on the torus for every $\varepsilon$.
We  assume some regularity of $f$ with respect to $\varepsilon$ such that the power series expansion in
$\varepsilon$ is valid.

We regard the set $\RR\times\RR^{+}$ of points $(\alpha,\varepsilon)$ as the parameter space.
Typically, when $\varepsilon=0$ and $\alpha_{0}\Omega$ is resonant (see the definition in Subsection \ref{resonant intrinsic}),
the torus decomposes into a family of  co-dimension one tori. Each orbit is dense on such a lower dimension tori, but not
in $\TT^{d}$.

We are interested in the dynamical behavior of $F_{\alpha, \varepsilon}$ for those parameters $\alpha$ close to
$\alpha_{0}$ (such that $\alpha_{0}\Omega$ is  resonant)
and small $\varepsilon$. To this end, we develop a resonant normal form for the foliation
preserving torus maps. Assuming some non-degeneracy condition, we show that the
invariant surface for the resonant normal form persists under  perturbation. See
Theorem \ref{phase locking theorem}.

The fact that the resonance
module of 
$F_{\alpha, \varepsilon}$
is one-dimensional (or the natural frequency is zero)
leads to the fact that the dynamics of such maps
is very different from that of generic maps of the torus
close to rotations. As we will see, when the natural frequency is
not zero, 
the structures that are generated by  resonances are codimension one
manifolds which obstruct the foliation.  For generic maps
this consequence is not true and one can get many more
objects generated by resonances.  As we will see, this
difference appears already in  the study of formal power series solutions.

Furthermore, due to preserving the $\Omega$-direction, $F_{\alpha,\varepsilon}$ can also
be described as a skew product flow in the neighborhood of the invariant surface generated by
resonance. Then, the evolution on the normal coordinate can be linearized under some
hyperbolic hypothesis. See also \cite{KP90}. Moreover, if the evolution on the invariant
circle can be conjugated to a rigid rotation, then the skew product flow can be further
reduced to linear system with constant coefficients.
See Theorem \ref{Sternberg} and Theorem \ref{constant expanding}.

Finally, we give a simple example to show our theoretical results.

\subsection{Invariant surface generated by resonance}
In the following, we assume that $\alpha_{0}\Omega$ is resonant. Then Proposition
\ref{one dimension resonance} implies that there exists a matrix
$\mathfrak{A}\in\textrm{SL}(d,\ZZ)$, $\omega\in\RR^{d-1}$ and $L\in\ZZ^{d}$ such that
\begin{equation}\label{intrinsic frequency}
 \mathfrak{A}\alpha_{0}\Omega=
\begin{pmatrix} \omega\\ 0\end{pmatrix}+ L
\end{equation}
with $\omega \cdot m\neq 0$ for any $m\in\ZZ^{d-1}\setminus\{0\}$.
Furthermore, we assume the intrinsic frequency $\omega$ satisfies Diophantine condition
\eqref{Diophantine condition}
\begin{equation*}
|k\cdot \omega-n|>\nu |k|^{-\tau}
\end{equation*}
for all $k\in \mathbb{Z}^{d-1}\setminus\{0\}$ and  $n\in \mathbb{Z}$.

\subsubsection{Resonant normal form}

We apply the averaging method to obtain a resonant normal form for \eqref{original system}
in two steps.
We first make some heuristic calculations to formulate the resonant normal form. Then we give
a detailed analysis on the convergence of Fourier series which requires solving the cohomology equation
\eqref{cohomology-torus-1}.

To this end, making change of variables by $H_{\varepsilon}=Id+ \varepsilon h_{\varepsilon}\Omega$, we have
\begin{equation}\label{conjugation}
\begin{aligned}
H_{\varepsilon}^{-1}&\circ F_{\alpha, \varepsilon}\circ H_{\varepsilon}-T_{\alpha\Omega}\\
=&
\varepsilon \left(h^{0}-h^{0}\circ T_{\alpha\Omega}+f^{0}\right)\Omega\\
&+\varepsilon^{2}\left(h^{1}-h^{1}\circ T_{\alpha\Omega}+f^{1}+Df^{0}
\cdot\Omega \cdot h^{0}
\right.\\
&\left.\quad+Dh^{0}\circ T_{\alpha\Omega}\cdot\Omega\cdot h^{0}\circ T_{\alpha\Omega}-Dh^{0}\circ T_{\alpha\Omega}\cdot\Omega\cdot (h^{0}+f^{0})
\right)\Omega\\
&+\cdots ,
\end{aligned}
\end{equation}
where $T_{\alpha\Omega}(x)=x+\alpha\Omega$ for $x\in\TT^{d}$,
$f_{\varepsilon}$ and $h_{\varepsilon}$ are expanded into power series in $\varepsilon$
as
$$
f_{\varepsilon}=f^{0}+\varepsilon f^{1}+\cdots,\quad
h_{\varepsilon}=h^{0}+\varepsilon h^{1}+\varepsilon^{2}h^{2}+\cdots.
$$
For those $\alpha$ close enough to $\alpha_{0}$, we denote
\begin{equation}\label{difference}
\begin{split}
 \Delta\alpha=&\alpha-\alpha_{0},\\
 \widetilde{\Delta}(\alpha)=&  H_{\varepsilon}^{-1}\circ F_{\alpha,\varepsilon}\circ H_{\varepsilon}
-H_{\varepsilon}^{-1}\circ F_{\alpha_{0},\varepsilon}\circ H_{\varepsilon}.
\end{split}
\end{equation}

From the facts that $F_{\alpha, \varepsilon}-F_{\alpha_{0}, \varepsilon}=\Omega\Delta\alpha$ and
$D^{n}(F_{\alpha,\varepsilon}\circ H_{\varepsilon})$ is independent of $\alpha$,
we have
\begin{equation}\label{difference estimate}
\|\widetilde{\Delta}(\alpha)\|_{\mathscr{C}^{0}}\leq \|D(H_{\varepsilon}^{-1})\|_{\mathscr{C}^{0}}
\cdot|\Omega|\cdot|\Delta\alpha|
\end{equation}
and
\begin{equation}\label{diff est 2}
\begin{split}
\|D^{n}(\widetilde{\Delta}(\alpha))\|_{\mathscr{C}^{0}}\leq &
C_{n} \|H_{\varepsilon}^{-1}\|_{\mathscr{C}^{n+1}}\cdot \|F_{\alpha,\varepsilon}\circ H_{\varepsilon}\|_{
\mathscr{C}^{n}}\\
& \times  (1+\|F_{\alpha,\varepsilon}\circ H_{\varepsilon}\|_{
\mathscr{C}^{n}}^{n-1})\cdot|\Omega|\cdot|\Delta\alpha|
\end{split}
\end{equation}
for any $n\geq 1$,
where $\|\cdot\|_{\mathscr{C}^{n}}$ is defined in \eqref{norm c n}.

The estimates in \eqref{difference estimate} and \eqref{diff est 2} are standard and can be found in \cite{dlLO99}
in a uniform way. Nevertheless, \eqref{difference estimate} shows that it suffices to find a normal form
for $F_{\alpha_{0},\varepsilon}$ when $\Delta\alpha$ is sufficiently small.

In order to simplify \eqref{conjugation} with $\alpha=\alpha_{0}$, we solve the following
\emph{cohomology equations} up to some resonant terms
\begin{equation}\label{cohomo equ}
\begin{split}
 h^0\circ T_{\alpha_{0}\Omega}-h^0=&f^{0},\\
 h^{1}\circ T_{\alpha_{0}\Omega}-h^{1}=&f^{1}+[Df^{0}\cdot h^{0}
+Dh^{0}\circ T_{\alpha_{0}\Omega}\cdot h^{0}\circ T_{\alpha_{0}\Omega}\\
&-Dh^{0}\circ T_{\alpha_{0}\Omega}\cdot (h^{0}+f^{0})] ~\Omega,\\
\cdots&
\end{split}
\end{equation}
It is worth noticing that the L.H.S. of equations in \eqref{cohomo equ} (for
$h^{0}, h^{1},\cdots $) has the same form
while the R.H.S. are known depending on lower order terms
which can be  obtained step by step.

When expanding both sides of
\eqref{cohomo equ} into Fourier series directly and comparing them, it turns out that, for those $k$'s not lying on
the resonance module, the $k$-th  Fourier coefficient of the unknowns $h^{j}$ can be solved explicitly.

To be precise, we make  further linear transformation on equation \eqref{cohomo equ} such that
the intrinsic frequency $\omega$ in \eqref{intrinsic fre} plays a role. When  composing with an appropriate linear transformation
$\mathfrak{A}^{-1}$ on
the right, the first equation in \eqref{cohomo equ} reads
\begin{equation*}
 h^{0}\circ T_{\alpha_{0}\Omega}\circ \mathfrak{A}^{-1}-h^{0}\circ \mathfrak{A}^{-1}=(h^{0}\circ
\mathfrak{A}^{-1})\circ T_{\mathfrak{A}\alpha_{0}\Omega}
-h^{0}\circ \mathfrak{A}^{-1}=f^{0}\circ \mathfrak{A}^{-1}.
\end{equation*}
Then expanding $g^{0}=h^{0}\circ \mathfrak{A}^{-1}$ and $p^{0}=f^{0}\circ \mathfrak{A}^{-1}$ into Fourier series,
 we have
\begin{equation}\label{g Four coeff}
 (e^{2\pi\mathrm{i}\langle k,(\omega,0)\rangle}-1)\cdot \widehat{g}^{0}(k)
=\widehat{p}^{0}(k),
\end{equation}
where $\widehat{(\cdot)}(k)$ denotes the $k$-th Fourier coefficient.

Let $k=(m,l)$ with $m\in\ZZ^{d-1}$ and $l\in\ZZ$. For $m\neq 0$, it follows that
\[
 \widehat{g}^{0}(k)=\frac{\widehat{p}^{0}(k)}{e^{2\pi\mathrm{i}\langle m, \omega\rangle}-1}.
\]
For $m=0$, we choose $\widehat{g}^{0}(k)=0$ and leave equation \eqref{g Four coeff} unsolved.
Combining the choice of  the Fourier coefficients $\widehat{g}^{0}(k)$, we see that
\begin{equation*}
  h^{0}\circ \mathfrak{A}^{-1}-h^{0}\circ T_{\alpha_{0}\Omega}\circ \mathfrak{A}^{-1}+f^{0}\circ
\mathfrak{A}^{-1}=\sum_{l\in\ZZ}\widehat{p}^{0}((0,l))~ e^{2\pi\mathrm{i}\langle (0,l), \cdot\rangle}
\end{equation*}
and
\[
 h^{0}=g^{0}\circ \mathfrak{A}=\sum_{k=(m,l)\in \mathbb{Z}^{d-1}\times \mathbb{Z}, ~m\neq 0}
\widehat{g}^{0}(k)~ e^{2\pi\mathrm{i}\langle k,\mathfrak{A}\cdot\rangle}.
\]

When $F_{\alpha,\varepsilon}$ is analytic such that $f^{0}\circ \mathfrak{A}^{-1}\in\mathscr{A}_{\rho}$,
the analyticity of $h^{0}$ follows from the Diophantine condition \eqref{Diophantine condition}
and the exponential decay estimate  \eqref{exponential decay} for the Fourier coefficient.
 More precisely, for any $x\in\CC^{d}$ with $|\mathrm{Im}x|< \rho-\delta$, we have
\begin{equation}\label{cohomo analytic}
\begin{split}
 |g^{0}(x)|=&\sum_{0\neq m\in\ZZ^{d-1}}\sum_{l\in\ZZ}\widehat{g}^{0}(k) e^{2\pi\mathrm{i}\langle (m,l),x\rangle}\\
\leq&\sum_{0\neq m\in\ZZ^{d-1}}\sum_{l\in\ZZ}\frac{|\widehat{p}^{0}(k)|}{|e^{2\pi\mathrm{i}\langle m, \omega\rangle}
-1|}e^{2\pi (|m|+|l|)(\rho-\delta) }\\
\leq& \|f^{0}\circ \mathfrak{A}^{-1}\|_{\rho}\sum_{0\neq m\in\ZZ^{d-1}}\frac{\pi}{2}\nu^{-1}|m|^{\tau}e^{-2\pi|m|\delta}\sum_{l\in\ZZ}e^{-2\pi \delta |l|}\\
\leq & \nu^{-1}C_{d, \tau}\delta^{-(\tau+d)}\ \|f^{0}\circ \mathfrak{A}^{-1}\|_{\rho},
\end{split}
\end{equation}
where $\|\cdot\|_{\rho}$ is given by \eqref{analytic norm}.

Similarly, when $F_{\alpha, \varepsilon}$ belongs to $\mathscr{C}^{n}$, we have, by \eqref{poly decay}, that
\begin{equation}\label{cohomo diff}
 \begin{split}
  \|g^{0}\|_{\mathscr{C}^{n-s}}\leq &(2\pi)^{n-s}\sum_{0\neq m\in\ZZ^{d-1}}\sum_{l\in\ZZ}(|m|+|l|)^{n-s}
\frac{|\widehat{p}^{0}((m,l))|}{|e^{2\pi i \langle m, \omega\rangle}-1|}\\
\leq &\nu^{-1} C_{n} \sum_{0\neq m\in\ZZ^{d-1}}\sum_{l\in\ZZ}
(|m|+|l|)^{-s}\cdot |m|^{\tau} \cdot \|f^{0}\circ \mathfrak{A}^{-1}\|_{\mathscr{C}^{n}}\\
\leq &\nu^{-1} C_{n, d} \sum_{\mu\geq 0} \mu^{\tau-s+d-1}\cdot  \|f^{0}\circ \mathfrak{A}^{-1}\|_{\mathscr{C}^{n}},
\end{split}
\end{equation}
where $\|\cdot\|_{\mathscr{C}^{n}}$ is defined in \eqref{norm c n}.
The sum in the R.H.S. converges provided $s>\tau+d$.


Likewise,  denoting
 $$p^{1}\circ \mathfrak{A}=f^{1}+Df^{0}\cdot \Omega\cdot  h^{0}
+(Dh^{0}\cdot\Omega\cdot h^{0})\circ T_{\alpha_{0}\Omega})-Dh^{0}\circ T_{\alpha_{0}\Omega}\cdot\Omega\cdot(h^{0}+f^{0}),$$
 we obtain
\begin{equation*}
  h^{1}\circ \mathfrak{A}^{-1}-h^{1}\circ T_{\alpha_{0}\Omega}\circ \mathfrak{A}^{-1}+p^{1}
=\sum_{l\in\ZZ}\widehat{p}^{1}((0,l))~ e^{2\pi\mathrm{i} \langle (0,l), \cdot\rangle},
\end{equation*}
where $h^{1}$ can be written explicitly as $h^{0}$.

By induction, one can conduct the same averaging process to $h^{j+1}$ in \eqref{cohomo equ} when $h^{0},\cdots,h^{j}$ are already known but losing
some regularities at each step.
Above all, under the averaging procedure and the linear transformation,
one actually obtain the following normal form for $F_{\Omega_{0},\varepsilon}$, say, to the order of $N$
\begin{equation}\label{normal form}
\begin{aligned}
 &\mathfrak{A}\circ (H_{\varepsilon}^{N})^{-1}\circ F_{\alpha_{0}\Omega,\varepsilon}\circ H_{\varepsilon}^{N}\circ \mathfrak{A}^{-1}
\begin{pmatrix} x_{1}\\ x_{2}\end{pmatrix}\\
=&\begin{pmatrix} x_{1}\\ x_{2}\end{pmatrix}
+\begin{pmatrix}\omega\\ 0\end{pmatrix}+ \alpha_{0}^{-1} \begin{pmatrix} \omega+L_{1}\\ L_{2}\end{pmatrix} \sum_{j=1}^{N}\varepsilon^{j}
\sum_{l\in\ZZ} \widehat{p}^{j}((0,l))~ e^{2\pi i\langle l, x_{2}\rangle}\\
&+\varepsilon^{N+1}r(x_{1},x_{2};\varepsilon),
\end{aligned}
\end{equation}
where $(x_{1},x_{2})\in\torus^{d-1}\times\torus$, $L=(L_{1}, L_{2})\in \mathbb{Z}^{d}$ and $\varepsilon^{N+1}r(x_{1},x_{2};\varepsilon)$
is the Taylor remainder of order $\varepsilon^{N+1}$. Thus $r(x_{1},x_{2};\varepsilon)$
is bounded for sufficient small $\varepsilon$.
Furthermore, the invertible conjugacy function $H^{N}$ takes the form of
$$
H^{N}_{\varepsilon}=Id+\varepsilon h^{0}\Omega+\varepsilon^{2} h^{1}\Omega
+\cdots \varepsilon^{N} h^{N-1}\Omega.
$$

For future application, we introduce the  following notations. Typically,
 we decompose $r$ into
\begin{equation}\label{r}
r(x_{1},x_{2};\varepsilon)=\begin{pmatrix}r_{1}(x_{1},x_{2}; \varepsilon)\\
 r_{2}(x_{1},x_{2}; \varepsilon)\end{pmatrix}\in\TT^{d-1}\times\TT.
\end{equation}

Let
\begin{equation}\label{n}
\begin{aligned}
 n=\min\Big\{
 s\geq 1: ~&\sum_{l\in\ZZ}\widehat{p}^{j}((0,l))~ e^{2\pi \textrm{i}~ \langle l,\cdot\rangle}\equiv 0, ~ \forall 1\leq j\leq s,\\
&
\sum_{l\in\ZZ}\widehat{p}^{s}((0,l))~ e^{2\pi \textrm{i}~ \langle l,\cdot\rangle}\not\equiv 0
\Big\}
\end{aligned}
\end{equation}
and we denote
\begin{equation}\label{beta eta}
\begin{pmatrix}\beta(x_{2}, \varepsilon)\\
 \eta(x_{2}, \varepsilon)\end{pmatrix}= \alpha_{0}^{-1} \begin{pmatrix} \omega+L_{1}\\ L_{2}\end{pmatrix}
\sum_{n\leq j\leq n+m-1}\varepsilon^{j-n}
\sum_{l\in\ZZ} \ \widehat{p}^{j}((0,l))~ e^{2\pi \textrm{i} \langle l, x_{2}\rangle}
\end{equation}
for some integer $m\geq 1$.
By an abuse of notation, we also denote $\mathfrak{A}\circ\widetilde{\Delta}(\alpha)\circ\mathfrak{A}^{-1}$
by $\widetilde{\Delta}(\alpha)$ and decompose it into
\begin{equation}\label{delta hat omega}
\widetilde{\Delta}(\alpha)=\begin{pmatrix}\widetilde{\Delta}_{1}(\alpha)\\
 \widetilde{\Delta}_{2}(\alpha)\end{pmatrix}\in\TT^{d-1}\times\TT.
\end{equation}

Taking $N=n+m-1$, we summarize the averaging results in the following theorem.

\begin{theorem}\label{lemma normal form}
Given $p, q\in\NN$ with $p>q\geq (n+m-1)(d+\tau)$.
Assume $\alpha_{0}\Omega$ is resonant and the intrinsic frequency $\omega$
satisfies Diophantine condition \eqref{Diophantine condition}.
Assume further that $F_{\alpha, \varepsilon}(x)$ is real analytic $($or $\mathscr{C}^{p}$$)$ in $x$ and $\mathscr{C}^{n+m}$ in $\varepsilon$.
 Then, for sufficiently small
$\varepsilon$, the  torus
map $F_{\alpha,\varepsilon}$ can be  conjugated by a real analytic $($or $\mathscr{C}^{p-q}$, respectively$)$ invertible
 function $H_{\varepsilon}^{n+m-1}$ and a linear transformation
$\mathfrak{A}\in\textrm{SL}(d,\ZZ)$
to the following resonant normal form $\widetilde{F}_{\alpha, \varepsilon}: \TT^{d-1}\times\TT
 \rightarrow \TT^{d-1}\times\TT $ defined by
\begin{equation}\label{normal form equ}
\widetilde{F}_{\alpha,\varepsilon}
\begin{pmatrix}x_{1}\\x_{2}\end{pmatrix}=
\begin{pmatrix}
x_{1}+\omega+\varepsilon^{n}\beta(x_{2},\varepsilon)+
\varepsilon^{n+m}r_{1}(x_{1},x_{2}; \varepsilon)
+\widetilde{\Delta}_{1}(\alpha)
\\
 x_{2}
+\varepsilon^{n}\eta(x_{2},\varepsilon)
+\varepsilon^{n+m}r_{2}(x_{1},x_{2}; \varepsilon)
+\widetilde{\Delta}_{2}(\alpha)
\end{pmatrix},
\end{equation}
where $n,m$ and  $\beta,\eta,r_{1},r_{2},\widetilde{\Delta}_{1}(\alpha),\widetilde{\Delta}_{2}(\alpha)$
are given in \eqref{r}-\eqref{delta hat omega}.

Furthermore,  $r_{1}$ and $r_{2}$ are real analytic $($or $\mathscr{C}^{p-q}$, respectively$)$ periodic functions defined on $\TT^{d}$,
$\beta$ and $\eta$ are real analytic periodic $($or $\mathscr{C}^{p-q}$, respectively$)$ functions on $\TT$,
and $\widetilde{\Delta}(\alpha)$ satisfies \eqref{difference estimate}.
\end{theorem}

\subsubsection{Persistence of invariant surface}

By an abuse of notation, we still
denote $F_{\alpha, \varepsilon}$ by the resonant normal form developed in Theorem \ref{lemma normal form}.
In this section, we assume that there exists $x_{2}^*\in\TT$  such that
\begin{equation}\tag{\textbf{H1}}\label{zero point assumption}
\eta(x_{2}^*,0)=0\quad\textrm{and}\quad D_{1}\eta(x_{2}^{*},0)\neq 0.
\end{equation}
By the implicit function theorem,
we obtain $x_{2}^{*}(\varepsilon)$ such that $\eta(x_{2}^{*}(\varepsilon),\varepsilon)\equiv \eta(x_{2}^{*},0)$
 and $x^{*}_{2}(0)=x_{2}^*$.
Then it follows that $\Gamma=\{(x_{1},x_{2}^{*}(\varepsilon))\ |\ x_{1}\in\TT^{d-1}\}$ is an invariant surface
for the torus map
$(x_{1},x_{2})\mapsto(x_{1}+\omega+\varepsilon^{n}\beta(x_{2},\varepsilon),x_{2}
+\varepsilon^{n}\eta(x_{2},\varepsilon))$
on which the motion is a rotation.
Without loss of generality, we assume $x_{2}^{*}(\varepsilon)\equiv 0$.

In what follows, we show the existence of invariant surface for the resonant normal
form developed in Theorem \ref{lemma normal form}, which is close to $\Gamma$.
In the perturbative setting, we see that the invariant manifold
of $F_{\alpha,\varepsilon}$ can be represented by a graph
of $w: \TT^{d-1}\rightarrow\RR$, whose invariance determines a functional $\mathscr{F}$ defined
on some function space. Then by the contraction mapping arguments, we can  prove the existence
of a fixed point for $\mathscr{F}$, which corresponds to the desired invariant manifold of
 $F_{\alpha,\varepsilon}.$

\begin{theorem}\label{phase locking theorem}
 Let $F_{\alpha,\varepsilon}$ be a family of real analytic ($\mathscr{C}^{\infty}$
or finitely differentiable) torus maps given by \eqref{original system}.
Assume $\alpha_{0}\Omega$ is resonant with intrinsic frequency $\omega$ satisfying
Diophantine condition \eqref{Diophantine condition}.

Then if
the resonant normal form \eqref{normal form equ} satisfies non-degeneracy condition \eqref{zero point assumption},
 there exists a $(d-1)$-dimensional  finitely differentiable invariant torus
of $F_{\alpha, \varepsilon}$ for  those parameters $(\alpha,\varepsilon)$ in some neighborhood
of $(\alpha_{0},0)$.

If the map
is finitely  differentiable (but with enough regularity)
we have an analogue result
but the regularities of the low dimensional invariant torus may be
less  than those of the map.
\end{theorem}

In Appendix \ref{general torus map}, we prove a theorem on the existence of
invariant surface for a family of generic torus maps generated by the resonances, from which
Theorem \ref{phase locking theorem} is an immediate result. Moreover,
in the proof of Theorem \ref{phase locking theorem}, we will employ a particular
contraction mapping theorem (see \cite[Lemma 2.4]{HdlL16}),
which guarantees   the existence of the fixed point of the operator $\mathscr{F}$ by
verifying that $\mathscr{F}$ maps a closed subset in the high regularity space to itself, and
is a contraction in the low regularity space.
Similar ideas also appear in \cite{Lan73}.

\begin{remark}
In Theorem \ref{phase locking theorem}, we only obtain a finitely differentiable invariant
surface even though the map $F_{\alpha, \varepsilon}$ is $\mathscr{C}^{\infty}$ or analytic.

This depends a lot on the dynamics of the map in the invariant manifold. 
For instance, in the case that the manifold
is a one-dimensional circle and the dynamics
on the circle have an attractive periodic point,
it is shown in \cite{dlL97} that the 
invariant circles may be finitely differentiable
and the regularity is determined by the exponents of
the derivative at the periodic points. 
On the other hand, in the case that the maps inside of
the invariant circle have a Diophantine rotation number,
the circles are analytic when the map is analytic. 

In the case of higher dimensional invariant manifolds,
there are similar obstructions to regularity depending on what is
the dynamics inside of the manifold (it could have invariant manifolds
of further codimension).  This requires further study.
\end{remark}

\subsection{Computation of Lindstedt series.}
We formulate the Lindstedt series for the invariant surface of the foliation
preserving torus map generated by
resonance, which provides a formal power series of the invariant surface.

Consider the following foliation-preserving torus map
\begin{equation}\label{simple model}
F_{\varepsilon}(\theta)=\theta+\alpha_{0}\Omega+\varepsilon g(\theta)\Omega,
\end{equation}
where $\Omega$ is irrational and $\varepsilon>0$.
Assume also $\alpha_{0}\Omega$ is resonant.
Then from \eqref{intrinsic frequency}
we have
\[
\mathfrak{A}\circ F_{\varepsilon}\circ \mathfrak{A}^{-1}(\theta)=\theta+\left[
\begin{pmatrix} \omega\\ 0\end{pmatrix}+L\right]\cdot \left[1+\varepsilon \alpha^{-1}g(\mathfrak{A}^{-1}\theta)\right],
\]
and by an abuse of notation, we can write
\begin{equation*}
F_{\varepsilon}\begin{pmatrix} x\\ y\end{pmatrix}=\begin{pmatrix} x\\ y\end{pmatrix}+\begin{pmatrix} \omega\\ 0\end{pmatrix}
+\begin{pmatrix}\varepsilon(\omega+L_{1})g(x,y)\\ \varepsilon L_{2}g(x,y)\end{pmatrix}.
\end{equation*}

We look for $l_{\varepsilon}=(l_{\varepsilon}^{x}, l_{\varepsilon}^{y}): \TT^{d-1}\rightarrow \TT^{d}$ and $u_{\varepsilon}: \TT^{d-1}\rightarrow \TT^{d-1}$  such that
\begin{equation}\label{embedding}
F_{\varepsilon}\circ l_{\varepsilon}=l_{\varepsilon}\circ u_{\varepsilon},
\end{equation}
by which the graph of $l_{\varepsilon}$ is the desired invariant surface.
To formulate the Lindstedt series, we
expand $l_{\varepsilon}$ and  $u_{\varepsilon}$ into power series in $\varepsilon$ as
\[
l_{\varepsilon}(\sigma)=\sum_{j=0}^{\infty} l_{j}(\sigma) \varepsilon^{j}=\sum_{j\geq 0}\begin{pmatrix}l^{x}_{j}(\sigma)\\ l_{j}^{y}(\sigma)
 \end{pmatrix}
\varepsilon^{j},\quad
u_{\varepsilon}(\sigma)=\sum_{j\geq 0} u_{j}(\sigma)\varepsilon^{j},
\]
and then take formal calculations.

By matching the coefficients of the $\varepsilon^{0}$-terms in \eqref{embedding}, we obtain
$F_{0}\circ l_{0}=l_{0}\circ u_{0}$ or equivalently
\begin{equation*}
 \left\{
\begin{aligned}
&l_{0}^{x}(\sigma)+\omega=l_{0}^{x}(u_{0}(\sigma)),\\
&l_{0}^{y}(\sigma)=l_{0}^{y}(u_{0}(\sigma)),
\end{aligned}
\right.
\end{equation*}
which can be solved by choosing
\[
 l_{0}^{x}(\sigma)=\sigma, \quad l_{0}^{y}(\sigma)=y_{0}=\textrm{Constant},\quad u_{0}(\sigma)=\sigma+\omega.
\]
For those $\varepsilon^{1}$-terms, we get
\[
F_{1}\circ l_{0}(\sigma)+DF_{0}\circ l_{0}(\sigma)\cdot l_{1}(\sigma)=l_{1}\circ u_{0}(\sigma)+Dl_{0}\circ u_{0}(\sigma)\cdot u_{1}(\sigma),
\]
or equivalently
\begin{equation}\label{Lindstedt 1st}
 \left\{
\begin{aligned}
&l_{1}^{x}(\sigma+\omega)-l_{1}^{x}(\sigma)=(\omega+L_{1})g(\sigma,y_{0})-u_{1}(\sigma),\\
&l_{1}^{y}(\sigma+\omega)-l_{1}^{y}(\sigma)=L_{2}g(\sigma,y_{0}).
\end{aligned}
\right.
\end{equation}
It suffices to choose $y_{0}$ such that
\begin{equation}\label{average vanishing}
 \int_{\TT^{d-1}} g(\sigma, y_{0})\ \textrm{d}\sigma=0,
\end{equation}
and
$u_{1}(\sigma)\equiv 0.$  Then $l_{1}=(l_{1}^{x},l_{1}^{y})$ can be solved from
\eqref{Lindstedt 1st} by using Lemma \ref{cohomology-torus}, but with
the average $$\langle l_{1}\rangle=\int_{\TT^{d-1}} l_{1}(\sigma)~ \mathrm{d} \sigma$$
to be specified.

To clarify the induction, we proceed to compute the equation for those
$\varepsilon^{2}$-terms, which reads
\[
F_{2}\circ l_{0}+DF_{1}\circ l_{0}\cdot l_{1}+DF_{0}\circ l_{0}\cdot l_{2}=
 l_{2}\circ u_{0}+Dl_{1}\circ u_{0}\cdot u_{1}+Dl_{0}\circ u_{0}\cdot u_{2},
\]
or equivalently
\begin{equation}\label{Lindstedt 2nd}
 l_{2}(\sigma+\omega)-l_{2}(\sigma)=
\begin{pmatrix}
 (\omega+L_{1})Dg(\sigma, y_{0})l_{1}(\sigma)-u_{2}(\sigma)\\
L_{2}Dg(\sigma, y_{0})l_{1}(\sigma)
\end{pmatrix}
.
\end{equation}
Let
 $l_{1}(\sigma)=\langle l_{1} \rangle +\widetilde{l}_{1}(\sigma)$,
in which $\widetilde{l}_{1}(\sigma)$ is uniquely determined by \eqref{Lindstedt 1st}.
Then the R.H.S. of
\eqref{Lindstedt 2nd} reads
\[
 \begin{pmatrix}
  \omega+L_{1}\\ L_{2}
 \end{pmatrix}
\cdot Dg(\sigma,y_{0})\widetilde{l}_{1}(\sigma)+
\begin{pmatrix}
 \omega+L_1\\ L_{2}
\end{pmatrix}
Dg(\sigma, y_{0}) \langle l_{1}\rangle-
\begin{pmatrix}
 u_{2}(\sigma)\\ 0
\end{pmatrix}.
\]
In order the average of  R.H.S of \eqref{Lindstedt 2nd} to vanish, we need to choose
parameters $\langle l_{1}\rangle=(\langle l_{1}^{x}\rangle, \langle l_{2}^{y}\rangle)$
and $u_{2}(\sigma)$. More precisely, we have
\begin{equation*}
 \begin{aligned}
&\langle D_{2}g(\cdot, y_{0})\rangle\cdot \langle l_{1}^{y}\rangle= \langle
Dg(\cdot, y_{0})\widetilde{l}_{1}(\cdot) \rangle,\\
 &(\omega+L_{1})\langle D_{2}g(\cdot, y_{0})\rangle\cdot \langle l_{1}^{y}\rangle-\langle u_{2}\rangle=
(\omega+L_{1})\langle Dg(\cdot, y_{0})\widetilde{l}_{1}(\cdot) \rangle ,\\
 \end{aligned}
\end{equation*}
since $\langle D_{1}g(\cdot, y_{0})\rangle=0$ by \eqref{average vanishing}.
Then if
\begin{equation}\label{average non-degeneracy}
 \langle D_{2}g(\cdot, y_{0})\rangle\neq 0,
\end{equation}
we just choose $u_{2}\equiv 0$ and take
\begin{equation*}
\langle l_{1}^{y}\rangle=\langle  D_{2}g(\cdot, y_{0})\rangle^{-1}\cdot
\langle
Dg(\cdot, y_{0})\widetilde{l}_{1}(\cdot) \rangle.
\end{equation*}

By induction, we assume that, under assumption \eqref{average vanishing} and \eqref{average non-degeneracy},
we can choose $u_{j}\equiv 0$ for all $1\leq j\leq n-1$ and always find
$\langle l_{j}^{y}\rangle$ for $1\leq j\leq n-2$ such that \eqref{embedding} holds up to $O(\varepsilon^{n})$.
Comparing the coefficients of $\varepsilon^{n}$ in
\eqref{embedding}, we obtain
\begin{equation*}
  \begin{split}
&\langle D_{2}g(\cdot, y_{0})\rangle\cdot \langle l_{n-1}^{y}\rangle= \langle
 \mathcal{G}_{n}[l_{0},\cdots, l_{n-1}; g] \rangle,\\
  &(\omega+L_{1})\langle D_{2}g(\cdot, y_{0})\rangle\cdot \langle l_{n-1}^{y}\rangle-\langle u_{n}(\sigma)\rangle=
(\omega+L_{1})\langle \mathcal{G}_{n}[l_{0},\cdots, l_{n-1}; g] \rangle,
 \end{split}
\end{equation*}
where $\mathcal{G}_{n}$ can be computed explicitly by the known
functions $g$ and $l_{0},\cdots, l_{n-1}$ from the induction procedure.
Then it suffices to choose $u_{n}(\sigma)\equiv 0$ and solve $\langle l_{n-1}^{y}\rangle$ from
the second equation above.
This completes the induction.

We conclude the above results in the following proposition.
\begin{proposition}\label{Lind}
 Let $\alpha_{0}\Omega$ be resonant with an intrinsic frequency $\omega\in \RR^{d-1}$ defined
 in \eqref{intrinsic frequency}. Assume there
exists $y_{0}\in\RR$ such that \eqref{average vanishing} and \eqref{average non-degeneracy} hold.
Then we can find a formal power series $l_{\varepsilon}=\sum_{j=01}^{\infty} l_{j} \varepsilon^{j}$ in $\varepsilon$ such that
\[
 F_{\varepsilon}\circ l_{\varepsilon}(\sigma)=l_{\varepsilon}(\sigma+\omega) .
\]
\end{proposition}

Combining Theorem \ref{phase locking theorem} and Proposition \ref{Lind}, we have
the following result for the foliation preserving torus map $F_{\varepsilon}$ in \eqref{simple model}.
\begin{theorem}\label{Lind theorem}
Let $p, q, N\in \mathbb{N}$ with $p-1> q\geq N(d+\tau)$.
Suppose $F_{\varepsilon}(\theta)$ is $\mathscr{C}^{p}$ in $\theta$ and is smooth in $\varepsilon$.
Under the assumptions of Proposition \ref{Lind}, there exist two $\mathscr{C}^{p-q-1}$ maps $l_{*}: \torus^{d-1}\rightarrow \torus$
 and $u_{*}: \torus^{d-1}\rightarrow \torus^{d-1}$ such that 
 \begin{equation*}
 F_{\varepsilon}\circ l_{*}= l_{*}\circ u_{*}.
 \end{equation*}
Let $l^{\leq N}=\sum_{j=0}^{N} l_{j}(\sigma) \varepsilon^{j}$  and $u_{0}(\sigma)=\sigma+\omega$
be the truncated Lindstedt series  obtained 
in Proposition \ref{Lind} satisfying
\begin{equation*}
\|F_{\varepsilon}\circ l^{\leq N}-l^{\leq N}\circ u_{0}\|_{\mathscr{C}^{p-q}}< C \varepsilon^{N+1}.
\end{equation*}
Then we have
\begin{equation*}
\|l^{\leq N}-l_{*}\|_{\mathscr{C}^{p-q-1}}\leq C' \varepsilon^{N+1},\quad
\|u_{0}-u_{*}\|_{\mathscr{C}^{p-q-1}}\leq C' \varepsilon^{N+1},
\end{equation*}
where the  constant $C'>0$ depends on $p, q, N, d, \tau $ and $g$.
\end{theorem}

\noindent\textbf{Proof.}~Since the conditions \eqref{average vanishing} and \eqref{average non-degeneracy}
verify \textbf{(H1)}, we obtain the invariant torus $\Gamma_{\varepsilon}$ of $F_{\varepsilon}$, which is
constructed by the graph of the map $l_{*}: \torus^{d-1}\rightarrow \torus$ in the proof
of Theorem \ref{phase locking theorem}. We denote by $u_{*}: \torus^{d-1}
\rightarrow \torus^{d-1}$ 
the evolution on the invariant surface $\Gamma_{\varepsilon}$. Moreover,
we have  that $l_{*}$ and $u_{*}$ are $\mathscr{C}^{p-q-1}$ smooth, in contrast with
the $\mathscr{C}^{p-q}$ resonant normal form of $F_{\varepsilon}$.
The remaining estimates follow directly from the construction of the Lindstedt series.
\qed

\begin{remark}
For a
 family of generic torus maps
 \begin{equation*}
 F_{\Omega, \varepsilon}(x)=x+\Omega+\varepsilon f(x),
\end{equation*}
 the Lindstedt series for the resonant invariant surface is more involved
 than the Lindstedt series we obtained for foliation preserving maps.
It is shown in Proposition \ref{appen Lind} in the Appendix that
 there exists $l_{\varepsilon}$ such that
\begin{equation*}
  F_{\Omega,\varepsilon}\circ l_{\varepsilon} (\sigma)= l_{\varepsilon}(
  \sigma+\omega+\varepsilon u_{1}+\varepsilon^{2} u_{2}+\cdots),
\end{equation*}
where  $\omega$ is the intrinsic frequency of
$\Omega$ and $\{u_{j}\}_{j\geq 1}$
is a sequence of constant vectors independent of $\sigma$.

Note that in the case of general maps of the torus, the frequency on the torus
depends on the perturbation. However, as we shall see in Subsection \ref{sect stru},
a $C^{1}$ perturbation of the foliation preserving torus map does not change the
rotational number of its uniformly attracting (or repelling) invariant circle.
\end{remark}

\subsection{Sternberg linearization around the invariant circle}\label{sect sternberg}
In this section,  we consider the particular case of $d=2$ for the analytic
foliation preserving torus map \eqref{original system} near resonance.
Then under the assumption of Theorem \ref{phase locking theorem}, there exists
an  invariant circle for all $(\alpha,\varepsilon)$ close to
$(\alpha_{0},0)$.

Now we study the local dynamics of $F_{\alpha,\varepsilon}$ around the invariant circle.
Up to the coordinate transformation, the restriction of $F_{\alpha,\varepsilon}$ on a small
neighborhood of the
invariant circle
can be characterized by a skew product map
\begin{equation}\label{skew prod}
  \varphi: \torus\times \RR\ni (\sigma,\rho)\mapsto (u(\sigma), \Gamma_{\sigma}(\rho))\in
  \TT\times \RR,
\end{equation}
where $\rho$ describes the normal coordinate and $\Gamma_{\sigma}(0)=0$.
As a result, the invariant circle in the $(\sigma,\rho)$-coordinate is characterized
by $\TT\times \{\rho=0\}$.

We recall that it could happen that the invariant circle is
finitely differentiable  only. On the other hand, due to
the preservation of the foliation, we have a set of normal 
fibers that are sent into each other by the map. The motion from
a fiber to its image is analytic. In mathematical language
this is described as a \emph{skew-product map}. The dynamics in
the base may be finitely differentiable but the dynamics on
the fiber is analytic.  Hence, for a description of
the result, it is important to consider functions
that have different regularity along different directions.

Assume the dynamics $u$ on the invariant circle is invertible and denote
\begin{equation}\label{w}
 w(\sigma)=u^{-1}(\sigma).
\end{equation}
Moreover, we assume the invariant circle is uniformly repelling. More precisely,
 assume that  in the normal coordinate,
there exists a constant  $\lambda>1$  such that for sufficiently small
 $\delta>0$, we can always find a $\gamma>0$
 such that
\begin{equation}\label{pinching condition}
|\Gamma'_{\sigma}(\rho)-\lambda|<\delta
\end{equation}
for all $\rho\in B_{\gamma}^{\CC}(0)$ and $\sigma\in\TT$, where
$B_{\gamma}^{\CC}(0)$ is the closed ball centered at zero with radius $\gamma$ in the complex plane.

\begin{theorem}{$($Sternberg Linearization Theorem$)$}\label{Sternberg}
Consider  the skew product map  $\varphi(\sigma,\rho)=(u(\sigma),\Gamma_{\sigma}(\rho))$
defined in \eqref{skew prod}. Assume $\Gamma_{\sigma}(\rho)$ is analytic in $\rho$ and $\Gamma_{\sigma}(0)=0$.
 Furthermore, $u(\sigma)$  and $\Gamma_{\sigma}(\rho)$ are continuous in $\sigma$
satisfying \eqref{w} and \eqref{pinching condition}.
Then there exists an invertible coordinate transformation $H: \TT\times\RR\rightarrow \TT\times\RR$ such that
\[
H^{-1}\circ \varphi\circ H(\sigma,\rho)=(u(\sigma), A_{\sigma}\rho),
\]
where $A_{\sigma}=\Gamma'_{\sigma}(0)$. Moreover, the transformation $H$ is analytic in $\rho$ and continuous in $\sigma$.
\end{theorem}

\noindent\textbf{Proof.} We denote
\begin{equation*}
h^{N}_{\sigma}=A_{\sigma}\cdot A_{w(\sigma)}\cdots A_{w^{N}(\sigma)}\circ \Gamma^{-1}_{w^{N}(\sigma)}\circ \cdots\circ \Gamma^{-1}_{\sigma},
\end{equation*}
and it follows that
\begin{equation*}
h^{N}_{w(\sigma)}\circ \Gamma_{\sigma}^{-1}=A_{\sigma}^{-1}\cdot h^{N+1}_{\sigma}.
\end{equation*}

If $\lim_{N\rightarrow \infty}h^{N}_{\sigma}$ exists and denote
\[
h_{\sigma}=\lim_{N\rightarrow \infty}h^{N}_{\sigma},
\]
we then have
\[
h_{w(\sigma)}\circ \Gamma_{\sigma}^{-1}=A_{\sigma}^{-1}\cdot h_{\sigma}
\]
and equivalently
\[
h_{\sigma}\circ \Gamma_{\sigma}\circ h^{-1}_{w(\sigma)}=A_{\sigma}.
\]
Let
\[
H: \TT\times\RR\ni (\sigma, \rho)\mapsto (\sigma, h_{w(\sigma)}^{-1}(\rho))\in\TT\times\RR.
\]
We immediately have
\[
H^{-1}\circ \varphi\circ H(\rho, \sigma)= (u(\sigma), A_{\sigma}\rho).
\]

The only thing left is to show the existence and analyticity of the
limit. Observing the fact that $A_{\sigma}$ and
$\Gamma_{\sigma}(\rho)$ are tangent, we have:

\begin{equation*}
\begin{aligned}
& \|h_{\sigma}^{N+1}-h_{\sigma}^{N}\|_{B_{\gamma}^{\CC}(0)}\\
=&
\left\|A_{\sigma}\cdot A_{w(\sigma)}\cdots A_{w^{N}(\sigma)}\left[
A_{w^{N+1}(\sigma)}\circ \Gamma^{-1}_{w^{N+1}(\sigma)}-Id
\right]\circ \Gamma^{-1}_{w^{N}(\sigma)}\circ \cdots\circ \Gamma^{-1}_{\sigma}\right\|_{B_{\gamma}^{\CC}(0)}\\
\leq & (\lambda+\delta)^{N}\cdot C ~\|\Gamma^{-1}_{w^{N}(\sigma)}\circ \cdots\circ \Gamma^{-1}_{\sigma}\|^{2}_{B_{\gamma}^{\CC}(0)}.
\end{aligned}
\end{equation*}
Moreover, since $\Gamma_{\sigma}(0)=0$, the mean value theorem implies
\[
\|h_{\sigma}^{N+1}-h_{\sigma}^{N}\|_{B_{\gamma}^{\CC}(0)}\leq C\left[\dfrac{\lambda+\delta}{(\lambda-\delta)^{2}}\right]^{N}.
\]
Therefore, the sequence $\{h^{N}_{\sigma}\}$ converges uniformly to an analytic function $h_{\sigma}(\rho)$ (analytic in $\rho$) provided
$\delta$ is small such that
\[
\lambda+\delta<(\lambda-\delta)^{2}.
\]
This completes the proof of Theorem \ref{Sternberg}.
\qed

In the case that the dynamics $u$ on the invariant circle can be conjugated to a Diophantine
rotation, the non-autonomous Sternberg Linearization Theorem \ref{Sternberg} can also be strengthened in such
a way that, in the normal coordinate, the effect of $\varphi$ is simply to expand by a
constant factor.

More precisely, we give a general result Theorem \ref{constant expanding} which may be of interest itself.
This theorem tells that, when the map in the base is conjugate
to a Diophantine rotation we can get the normal contraction
to be constant  (by a smooth change of variables). 

\begin{theorem}\label{constant expanding}
Let $\psi: \TT^{d-1}\times\RR\rightarrow \TT^{d-1}\times\RR$ be of the form
\begin{equation*}
 \psi(\sigma,t)=(v(\sigma), a(\sigma)t )
\end{equation*}
with $v\in\scr{C}^{r}(\TT^{d-1},\TT^{d-1})$ a diffeomorphism and the function $a(\sigma)$ belonging
to $\scr{C}^{r}(\TT^{d-1}, \RR\setminus\{0\})$.
Assume that there exists a diffeomorphism $h\in\scr{C}^{r}(\TT^{d-1},\TT^{d-1})$ and a vector
$\Omega\in\scr{D}_{d-1}(\nu, \tau)$ such that
\begin{equation}\label{conjugacy u}
 h^{-1}\circ v\circ h(\sigma)=\sigma+\Omega.
\end{equation}
Assume further that $\tau<r$.

Then, there exists a function $b\in\scr{C}^{r-\tau}(\TT^{d-1}, \RR\setminus\{0\})$ and a constant
$\kappa\in\RR\setminus\{0\}$ such that the diffeomorphism
\begin{equation*}
 H(\sigma, t)=( h(\sigma), b(\sigma)t)
\end{equation*}
satisfies
\begin{equation}\label{conjugacy constant}
 \psi\circ H(\sigma, t)=H(\sigma+\Omega, \kappa t).
\end{equation}

\end{theorem}

\noindent\textbf{Proof.}
Expanding \eqref{conjugacy constant}, we see that the first component of \eqref{conjugacy constant}
is just our assumption \eqref{conjugacy u}.
The second component of \eqref{conjugacy constant} is
\begin{equation*}
 a\circ h(\sigma)\cdot b(\sigma)= b(\sigma+\Omega)\cdot \kappa .
\end{equation*}

Noticing that $a(\sigma)\neq 0$, hence, it has the same sign for all $\sigma\in\TT^{d-1}$.
In the case that $a$ is always positive, we see that
\begin{equation}\label{cohomology ln b}
 (\log b)(\sigma+\Omega)-(\log b)(\sigma)=\log(a\circ h)(\sigma)-\log \kappa,
\end{equation}
which is a very standard cohomology equation in KAM theory for $\log b$ (see \cite[Lemma 2.1.8]{dlL01}).

Then we just take
\[
 \kappa=\exp\left(\int_{\TT^{d-1}}\log (a\circ h)(\sigma)\ \textrm{d}\sigma \right)
\]
and solve \eqref{cohomology ln b} for $\log b$ (see \cite{Rus75} for estimates and obtain
$b$ by exponentiating this result).

In the case that $a(\sigma)$ is negative, we have
\[
 (\log b)(\sigma+\Omega)-(\log b)(\sigma)=[\log (-a\circ h)](\sigma)-\log(-\kappa),
\]
and once again, we can repeat the analysis above.
This completes the proof of Theorem \ref{constant expanding}.
\qed

The meaning of  \eqref{conjugacy constant} in this paper is that, performing change of variable after
Theorem \ref{Sternberg}, we reduce the mapping $\varphi$ in the neighborhood of invariant
circle to the constant coefficient map $(\sigma, \rho)\rightarrow (\sigma+\Omega, \kappa \rho)$.

\begin{remark}\label{abundance of conjugation}
 It is worth noticing that if we consider a family of problems, we expect that
the hypothesis of existence of function $h$ conjugating $v$ to a rotation will be
satisfied by a positive measure set of parameters.
\end{remark}

\begin{remark}
 For general maps of the torus, it is expected that in the vicinities of higher multiplicity
resonance, one could find higher co-dimension torus and develop analogues of Sternberg
theorems \cite{KP90}.
 In the case of foliation preserving torus maps, as indicated in \cite{SdlL12}, only resonance of
multiplicity one happen.
\end{remark}

\subsection{Structural stability}\label{sect stru}

In Subsection \ref{sect sternberg}, we show the reducibility of the
skew product flow \eqref{skew prod} around the repulsive invariant circle to a linear flow
provided that the rotation number of the evolution on the invariant circle
is Diophantine. In this part, we study the structural stability of
the foliation preserving torus maps in $\torus^{2}$.

Let $T_{f}: \torus^{2}\rightarrow \torus^{2}$ be a
$C^{1}$ map of form \eqref{foliationpreserving} with $f>0$.
Assume that $T_{f}$ has a uniformly
attracting invariant circle $\Gamma_{f}\subset \torus^{2}$.
Recall that the only direction with possible non-vanishing Lyapunov exponent is
the direction of $\omega$. Hence, the invariant circle $\Gamma_{f}$  is transversal
to $\omega$. Moreover,
since $T_{f}$ has
no periodic points, we see that the rotation number
of $T_{f}|_{\Gamma}$ has to be irrational.

Let $g: \torus^{2}\rightarrow \real$ be  sufficiently
$C^{1}$ close to $f$ and remain positive. Consider the map $T_{g}$ given by
\eqref{foliationpreserving}, by the theory of normally hyperbolic
invariant manifold \cite{Fen73, Fen77}, it follows that $T_{g}$
also admits a uniformly attracting invariant circle $\Gamma_{g}$ close to $\Gamma_{f}$.

Since  $T_{g}|_{\Gamma_{g}}$ depends continuously on $g$, if the rotation number
of $T_{g}|_{\Gamma_{g}}$ is different from that of $T_{f}|_{\Gamma_{f}}$, we can
find a map $h$ in the homotopic functions
 $\{s f+(1-s)g: 0\leq s\leq 1\}$ connecting $f$ and $g$ such that
$T_{h}|_{\Gamma_{h}}$ has rational rotation number, which is impossible since
$T_{h}$ also has no periodic points.
This shows that small $C^{1}$ perturbation of the foliation
preserving torus map does not change the rotation number of its
uniformly attracting (or repelling) invariant circle.

\subsection{Example}\label{example}
In this part, we study a simple example of the foliation preserving torus map to
illustrate our theoretical results in this section.
Consider
\begin{equation}\label{example model}
F_{\varepsilon}(\theta)= \theta+\Omega+\varepsilon \Omega g(\theta),\quad \theta\in\torus^{2},
\end{equation}
where $\Omega=(\omega, 1)^{T}\in \real^{2}$, $\varepsilon>0$ and $g: \torus^{2}\rightarrow \real$ is given by
\begin{equation*}
g(\theta)=a+ \delta_{1} \sin (2\pi x)+\delta_{2} \sin(2\pi y),\quad \theta=(x,y).
\end{equation*}
Obviously, $\Omega$ is resonant and can be written as $\Omega=(\omega, 0)^{T}+(0, 1)^{T}$,
where $\omega$ is the intrinsic frequency of $\Omega$.
Then we have
\begin{equation*}
F_{\varepsilon}\begin{pmatrix} x\\ y\end{pmatrix}=\begin{pmatrix} x\\ y\end{pmatrix}+\begin{pmatrix} \omega\\ 0\end{pmatrix}
+\begin{pmatrix}\varepsilon \omega g(x,y)\\ \varepsilon g(x,y)\end{pmatrix}.
\end{equation*}

Assuming $\omega$ is Diophantine, we obtain from Theorem \ref{lemma normal form}  the resonant normal form
\begin{equation*}
F_{\varepsilon}(x, y)=\begin{pmatrix}
x+\omega+\varepsilon \beta(y) +\varepsilon^{2} r_{1}(x,y;\varepsilon)\\
y+\varepsilon \eta(y)+\varepsilon^{2} r_{2}(x,y;\varepsilon)
\end{pmatrix},
\end{equation*}
where
\begin{equation*}
\eta (y)=a+ \delta_{2} \sin(2\pi y), \quad \beta(y)= \omega (a+\delta_{2}\sin (2\pi y)).
\end{equation*}
If $|a|< |\delta_{2}|$, there exists a number
$y_{*}\in [0,1]$ such that
\begin{equation} \label{invariance_average} 
  \eta(y_{*})=0
  \end{equation}
and $\eta'(y_{*})=2\pi \delta_{2} \cos (2\pi y_{*})\neq 0$, which verifies the non-degeneracy condition
\textbf{(H1)} in Theorem \ref{phase locking theorem}.
It then follows that $F_{\varepsilon}$ has a one-dimensional invariant circle for sufficiently small
$\varepsilon$.

\begin{remark}\label{twosolutions}
  Note that, in general, we will get two numbers $y_*$ solving $\eta(y_*) = 0$:
  one of them with $\eta'(y_*) > 0$ and another one with
  $\eta'(y_*) < 0$.
We will continue the analysis for one of them but the calculations
will apply to the other one just as well.
\end{remark}

Next, we employ the Lindstedt series to study the invariant circle generated by the resonance.
Obviously, conditions \eqref{average vanishing} and \eqref{average non-degeneracy} are exactly the
non-degeneracy condition \textbf{(H1)}. From Proposition \ref{Lind}, we have
\begin{equation*}
F_{\varepsilon}\circ l_{\varepsilon}(\sigma)= l_{\varepsilon}(\sigma+\omega),
\end{equation*}
where
\begin{equation*}
l_{\varepsilon}(\sigma)=\begin{pmatrix}
l^{x}_{\varepsilon}(\sigma)\\ l^{y}_{\varepsilon}(\sigma)
\end{pmatrix}=
\begin{pmatrix}
\sigma+\varepsilon l^{x}_{1}(\sigma)+ \varepsilon^{2} l^{x}_{2}(\sigma)+\cdots\\
y_{*}+\varepsilon l^{y}_{1}(\sigma)+ \varepsilon^{2} l^{y}_{2}(\sigma)+\cdots
\end{pmatrix}.
\end{equation*}
It follows that the graph $\mathfrak{C}$ of $l_{\varepsilon}$ is the invariant circle  of $F_{\varepsilon}$,
whose existence has already been established. Introduce the local coordinate $(\rho, \sigma)$
around the invariant circle $\mathfrak{C}$ by  defining
\begin{equation*}
\Phi: \real\times \torus\ni (\rho, \sigma)\mapsto (l^{x}_{\varepsilon}(\sigma)+\omega \rho, l^{y}_{\varepsilon}(\sigma)+ \rho)=(x,y)\in \torus^{2}.
\end{equation*}
Obviously, $\det\frac{\partial (x, y)}{\partial (\rho, \sigma)}=-1+O(\varepsilon)\neq 0$.
Let
\begin{equation*}
\Phi^{-1}(x, y)= \begin{pmatrix} y-y_{*}\\ x-\omega(y-y_{*})
\end{pmatrix}+\varepsilon \begin{pmatrix} \varphi(x, y) \\ \psi(x, y)
\end{pmatrix}+ O(\varepsilon^{2}).
\end{equation*}
It follows from $\Phi^{-1}\circ \Phi=Id$ that
\begin{equation*}
\varphi(\sigma+\omega\rho, y_{*}+\rho)=l_{1}^{y}(\sigma),
\end{equation*}
which is independent of $\rho$. One easily solves from $\omega\cdot \partial_{x} \varphi+ \partial_{y} \varphi=0$ that
\begin{equation*}
\varphi(x, y)= l_{1}^{y}(x-\omega(y-y_{*})).
\end{equation*}
Considering $\tilde{\varphi}= \varphi\circ F_{\varepsilon}\circ \Phi$, we compute
\begin{equation*}
\frac{d\tilde{\varphi}}{d\rho}= (\omega+O(\varepsilon)) \partial_{x} \varphi+ (1+O(\varepsilon)) \partial_{y}\varphi=O(\varepsilon).
\end{equation*}
As a result, under the new coordinate, we have
\begin{equation*}
\Phi^{-1}\circ F_{\varepsilon}\circ \Phi (\rho, \sigma)= (\Gamma_{\sigma} (\rho), u(\sigma)),
\end{equation*}
where $\Gamma_{\sigma}(0)=0$ and
\begin{equation*}
\Gamma_{\sigma}(\rho)= \rho+l_{\varepsilon}^{y}(\sigma)+\varepsilon g\circ \Phi- y_{*}+\varepsilon \tilde{\varphi}(\rho,\sigma)
+O(\varepsilon^{2}).
\end{equation*}
Compute
\begin{equation*}
\frac{d}{d\rho}\Big {|}_{\rho=0} (g\circ \Phi)=2\pi \left[\delta_{1}\omega \cos(2\pi l_{\varepsilon}^{x}(\sigma))+\delta_{2}
\cos(2\pi l^{y}_{\varepsilon}(\sigma))\right].
\end{equation*}
Recalling that $l^{y}_{\varepsilon}(\sigma)= y_{*}+O(\varepsilon)$ and $\cos(2\pi y_{*})=\pm \sqrt{1-(a/\delta_{2})^{2}}$,
we can always choose $\delta_{1}$ and  $\delta_{2}$ with $2\omega |\delta_{1}|<|\delta_{2}| \sqrt{1-(a/\delta_{2})^{2}}=\lambda $ 
and $\delta_{2} \cos (2\pi y_{*})>0$ such that
$\frac{d}{d\rho}|_{\rho=0} (g\circ \Phi)=\lambda\pi+O(\varepsilon)\neq 0$.
As a result, for  $0<\varepsilon\ll 1$, we have
\begin{equation*}
\Gamma_{\sigma}'(0)> 1+\frac{\lambda\pi}{4} \varepsilon>1,
\end{equation*}
which implies the invariant circle $\mathfrak{C}$ is uniformly repelling.

As indicated in Remark~\ref{twosolutions}, there are
other solutions of the equation \eqref{invariance_average}.
In the preceding we have presented the calculations for
one of them with $\eta'(y_*) >0 $ -- this is the solution
that plays a more important role in the applications to delay equations
in \cite{dlL21b-TBA}. Proceeding in the same way with the other solution, we
obtain another invariant circle which is attractive. 

For a  family of foliation preserving torus maps
$
F_{\alpha,\varepsilon}(\theta)=\theta+\alpha \Omega+ \varepsilon\Omega g(\theta)
$
close to $F_{\varepsilon}$ given by \eqref{example model} (i.e., $|\alpha-1|\ll 1$),
there also exists a uniformly repelling invariant circle $\mathfrak{C}_{\alpha}$ close to $\mathfrak{C}$
by using Theorem \ref{phase locking theorem} when $0<\varepsilon\ll1$.

Since the rotational number of $\mathfrak{C}_{\alpha}$ is continuous in $\alpha$, we see that, for
for many points  $\alpha$ close to one, Theorem \ref{constant expanding} holds and the map $F_{\alpha, \varepsilon}$ in the neighborhood of $\mathfrak{C}_{\alpha}$
can be reduced to an autonomous linear map.

A more delicate version of the argument -- which we do
not carry out in detail -- shows that the set of
parameters where  the rotation number is Diophantine is of
positive measure.

\section{KAM theory for foliation preserving torus maps}\label{sec KAM}

In this section, we consider the conjugation problem of the foliation preserving torus maps to
the rigid rotation by the KAM techniques.
Instead of treating the nonlinearity as a perturbation, we study a family of
non-perturbation foliation preserving torus maps.

The KAM theorem we present is in an a-posteriori format. That is, given an approximate solution satisfying
some non-degeneracy condition, there is a true solution nearby.

More precisely, we study the following foliation preserving torus maps
\begin{equation}
F(x)=x+\Omega f(x)
\end{equation}
defined on $\TT^{d}$ with non-resonant $\Omega$.
Assume $\alpha\Omega$ satisfy the Diophantine condition \eqref{Diophantine condition} and
$f$ is real analytic.

Following the extra parameter method of \cite{Mos67}, we find
an extra parameter $\lambda$ and a foliation preserving map in such
a way that the mapping $F+\lambda \Omega$ is conjugate to
the rotation by $\alpha$.

\begin{remark}
  We recall how the result thus obtained translates into results for
  families.
  
From the a-posteriori format
of the theorem, it follows that $\lambda(\alpha)$ is
a Lipschitz function for $\alpha$ defined in the set of
parameters.  If we consider a family of parameters indexed by
$\varepsilon$,  we obtain $\lambda(\varepsilon, \alpha)$ is
a Lipschitz function.

If we study the equation $\lambda(\varepsilon, \alpha) = 0$
using the Lipschitz implicit function theorem (under
the assumption that $\partial_\varepsilon \lambda \ne 0$), we obtain
that for a positive set of parameters $\varepsilon$, the
map can be conjugated to a rotation. A numerical study of this
problem for the families of maps that appear in the
study of cavities with moving boundaries appears in \cite{PdlLV03}. 
\end{remark} 

More concretely,
we are looking for a real analytic periodic function $h$ and parameter $\lambda\in\RR$ in
such a way that
\begin{equation}
 F\circ H=H\circ T_{\alpha\Omega}+\lambda\Omega,
\end{equation}
where $H=Id+h\Omega$ and $T_{\alpha\Omega}(x)=x+\alpha\Omega $.

Denoting the functional $\mathscr{F}$ by
\begin{equation}
\begin{split}
 \scr{F}[h,\lambda]\cdot \Omega=&F\circ H-H\circ T_{\alpha\Omega}-\lambda\Omega\\
=& [h-h\circ T_{\alpha\Omega}-(\alpha+\lambda)+f\circ (Id+h\Omega)]\cdot\Omega,
\end{split}
\end{equation}
 the conjugation problem of the foliation preserving torus maps is transformed into
finding the solution of the functional equation
\begin{equation}\label{KAM FE}
 \scr{F}[h, \lambda]=0.
\end{equation}
Following the conventional notations in the KAM theory, we denote in the sequel
various constants by the letter $C$ with some subscripts indicating
the dependence on the given quantities, which might be different from line to line.
These constants could be made explicit from the context, but need not be.

\begin{theorem}\label{KAM Theorem}
Let $\alpha\Omega$ satisfy Diophantine condition \eqref{Diophantine condition} and
$f$ belong to the analytic function space $\scr{A}_{\rho}$.
Assume that there is an approximate solution $(h_{0},\lambda_{0})$ satisfying
\begin{enumerate}[(i)]
 \item $h_{0}\in\scr{A}_{\rho}$;
 \item $DH_{0}=I+Dh_{0}\cdot\Omega$ is invertible with $(DH_{0})^{-1}\in\mathscr{A}_{\rho}$;
 \item $\det \langle I+ Dh_{0}\cdot \Omega \rangle \neq 0$ with $\langle \cdot\rangle$ being the average of periodic functions.

\end{enumerate}

Denoting the initial error by $e_{0}$, i.e.,
\begin{equation}\label{initial error}
e_{0}=\scr{F}[h_{0}, \lambda_{0}],
\end{equation}
then if $\|e_{0}\|_{\rho}$ is sufficiently small, there is a true solution $(h, \lambda)$ of
\eqref{KAM FE} satisfying $h\in\mathscr{A}_{\rho/2}$ and
\[
 |\lambda-\lambda_{0}|,\ \|h-h_{0}\|_{\rho/2}\leq C(4/\rho)^{2\tau} \|e_{0}\|_{\rho},
\]
where the constant $C$ depends only on the given quantities $h_{0}, \Omega, d$ and
$\nu$ in \eqref{Diophantine condition}.

\end{theorem}

We apply the Nash-Moser method to prove Theorem \ref{KAM Theorem}. In subsection
\ref{Newton equ sect.}, we analyze the Newton equation for the functional equation
\eqref{KAM FE} in which the small divisor problem is overcome by the classical
cohomology equation \eqref{cohomology-torus-1}. Then we prove the convergence of
the Newton iteration in subsection \ref{convergence proof}.

\subsection{Newton equation}\label{Newton equ sect.}
The Newton equation of \eqref{KAM FE} is
\begin{equation}\label{Newton equ}
 \Delta h-\Delta h\circ T_{\alpha\Omega}-\Delta\lambda+Df\circ (Id+h\Omega)\cdot \Omega \Delta h=-e .
\end{equation}
Differentiating both sides of
\begin{equation*}
 \scr{F}[h, \lambda]=e,
\end{equation*}
we obtain
\begin{equation}\label{derivative error}
 Dh-Dh\circ T_{\alpha\Omega}+Df\circ (Id+h\Omega) (I+ Dh\cdot \Omega)=De .
\end{equation}

Denoting
\begin{equation}\label{transformed correction}
 \Delta V=(I+Dh\cdot \Omega)^{-1}\cdot \Omega\cdot \Delta h
\end{equation}
and substituting \eqref{derivative error}-\eqref{transformed correction} into \eqref{Newton equ},
we have
\begin{align*}
\Delta V\circ T_{\alpha\Omega}-\Delta V
 = & (DH)^{-1}\circ T_{\alpha\Omega} \cdot e\cdot \Omega
-(DH)^{-1}\circ T_{\alpha\Omega}\cdot \Delta\lambda\cdot \Omega\\
+& (DH)^{-1}\circ T_{\alpha\Omega}
\cdot \Omega\cdot De \cdot \Delta V .
\end{align*}

Since $De\cdot \Delta V$ is quadratic in error, we ignore it for the moment and consider
the modified Newton equation
\begin{equation}\label{modified Newton}
 \Delta V\circ T_{\alpha\Omega}-\Delta V= \left[(DH)^{-1}\circ T_{\alpha\Omega} \cdot e
-(DH)^{-1}\circ T_{\alpha\Omega}\cdot \Delta\lambda \right] \cdot \Omega .
\end{equation}

Then, by Lemma \ref{cohomology-torus}, we can solve \eqref{modified Newton} uniquely
(in the sense of vanishing average)
provided
\begin{equation}\label{KAM non-degeneracy}
 \det\ \langle I+ \Omega Dh \rangle \neq 0,
\end{equation}
where $\langle \cdot\rangle $ denotes the average of periodic functions. Furthermore, we also have
\begin{equation*}
 \|\Delta V\|_{\rho-\delta}\leq C_{\nu, d}\cdot \|(DH)^{-1}\|_{\rho}\cdot \delta^{-\tau}\cdot \|e\|_{\rho},
\end{equation*}
and
\begin{equation}\label{delta lambda}
 |\Delta\lambda|=\left|
\langle (DH)^{-1} \rangle^{-1}\cdot \langle (DH)^{-1}\circ R_{\alpha\Omega}\cdot e \rangle
\right|\leq C_{\langle DH\rangle}\cdot  \|e\|_{\rho} .
\end{equation}
Hence, it follows from \eqref{transformed correction} that
\begin{equation}\label{delta h}
 \|\Omega \Delta h\|_{\rho-\delta}\leq C_{\nu, d } \|DH\|_{\rho}\cdot
\|(DH)^{-1}\|_{\rho}\cdot \delta^{-\tau}\cdot \|e\|_{\rho}.
\end{equation}

Now we give the estimates for the new error
\begin{equation*}
 e^{+}=\scr{F}[h^{+}, \lambda^{+}]
\end{equation*}
with
\begin{equation*}
 h^{+}=h+\Delta h,\quad \lambda^{+}=\lambda+\Delta \lambda.
\end{equation*}

Since
\begin{equation*}
\begin{split}
 \scr{F}[h^{+}, \lambda^{+}]=&\scr{F}[h, \lambda]+D\scr{F}[h,\lambda](\Delta h, \Delta \lambda)\\
&+\int_{0}^{1} s \int_{0}^{1} D^{2}\scr{F}[h+st\Delta h, \lambda, st\Delta \lambda]
(\Delta h, \Delta \lambda)^{\otimes 2}\
\textrm{d} t\ \textrm{d} s\\
=& (DH)^{-1}\circ R_{\alpha\Omega}\cdot \Omega \cdot De\cdot \Delta V\\
&+ \int_{0}^{1} s \int_{0}^{1} D^{2} f\circ (Id+ \Omega h+\Omega st \Delta h)\cdot (\Omega \Delta h)
^{\otimes 2} \textrm{d} t\ \textrm{d} s,
\end{split}
\end{equation*}
we see that
\begin{equation*}
\begin{split}
 \|\scr{F}[h^{+}, \lambda^{+}]\|_{\rho-\delta}\leq & C_{\nu, d} \cdot \|(DH)^{-1}\|_{\rho}^{2}\cdot
|\Omega|\cdot \delta^{-(\tau+1)}\cdot \|e\|_{\rho}^{2}\\
&\phantom{A}+ C_{\nu, d}^{2} \|D^{2}f\|_{\rho} \cdot
\|DH\|_{\rho}^{2}\cdot  \|(DH)^{-1}\|_{\rho}^{2}\cdot \delta^{-2\tau} \cdot \|e\|_{\rho}^{2}\\
\leq & C_{\nu, d, h_{0}, \Omega} \cdot \delta^{-2\tau}\cdot\|e\|_{\rho}^{2}.
\end{split}
\end{equation*}

We conclude the above analysis in  the following iterative lemma.

\begin{lemma}\label{iteration lemma}
Given the approximate solution
$(h, \lambda)\in \mathscr{A}_{\rho}\times\RR$ with the error $e=\scr{F}[h, \lambda]$.
 Assume $DH$ is invertible and \eqref{KAM non-degeneracy} holds. Then, for any $0<\delta<\rho$, there exists
$(\Delta h, \Delta \lambda)\in\mathscr{A}_{\rho-\delta}\times\RR$ such that
\[
 \|e^{+}\|_{\rho-\delta}=\|\scr{F}[h^{+}, \lambda^{+}]\|_{\rho-\delta}\leq C_{\nu, d, h, \Omega} \cdot
\delta^{-2\tau}\cdot \|e\|_{\rho}^{2} ,
\]
where  $h^{+}=h+\Delta h$ and  $\lambda^{+}=\lambda+\Delta \lambda.$
\end{lemma}

Note that the openness of the invertibility and
non-degeneracy condition \eqref{KAM non-degeneracy} in Lemma \ref{iteration lemma}
enables us to iterate the Newton steps.

Now we  start  proving the convergence of the iteration sequences as well as verifying the
conditions in Lemma \ref{iteration lemma} at each step. Another key point to ensure
the convergence is to prove the uniform boundedness of $DH_{n}$ and $(DH_{n})^{-1}$ in
the common analyticity domain developed below.
All the arguments are very standard in the KAM theory.

\subsection{Proof of the convergence}\label{convergence proof}

From the standard techniques in KAM theory, we use
the subscript $n$ to denote the $n$-th step for
the Newton iterations. More precisely,
 we choose the loss of
the analyticity domain $\sigma_{n}$ as
$
\sigma_{n}=2^{-(n-1)}\sigma
$
and $\sigma=\rho/4$. Let
$\rho_{n+1}=\rho_{n}-\sigma_{n+1}$ and
$\rho_{0}=\rho$.
Inductively, we assume the errors
$
e_{n}=\mathscr{F}[h_{n},\lambda_{n}]
$
satisfy
$\|e_{n}\|_{\rho_{n}}\leq \varepsilon_{n}$,
where $h_{0}, \lambda_{0}$ are given in \eqref{initial error} and
$\varepsilon_{0}=\|e\|_{\rho}$. Noted that $h_{n}$ and $\lambda_{n}$ are
inductively defined by $h_{n}=h_{n-1}+\Delta h_{n-1}$
and
$\lambda_{n}=\lambda_{n-1}+\Delta \lambda_{n-1}$.
Furthermore, we also
assume that
\begin{equation}\label{ind error}
\varepsilon_{n}=C\sigma_{n}^{-2\tau}\varepsilon_{n-1}^{2}.
\end{equation}
Generally, if \eqref{ind error} holds for all $n$,
it is easy to show that $\varepsilon_{n}$ approaches
zero when $\varepsilon_{0}$ is small enough.
Indeed, denoting $\tilde{\varepsilon}_{n}
=C\sigma^{-2\tau}2^{2\tau(n+1)}\varepsilon_{n}$,
from \eqref{ind error} one has $\tilde{\varepsilon}_{n+1}
=\tilde{\varepsilon}_{n}^{2}$, which implies
\begin{equation}
\tilde\varepsilon_{n}=[C(1/\sigma)^{2\tau}\varepsilon_{0}]^{2^{n}}.
\end{equation}
Then if
\begin{equation}
C(1/\sigma)^{2\tau}\varepsilon_{0}<1,
\end{equation}
$\varepsilon_{n}$ obviously approaches zero and satisfies
$$
\sum_{n=1}^{\infty}\tilde\varepsilon_{n}\leq
\sum_{n=1}^{\infty} [C(1/\sigma)^{2\tau}\varepsilon_{0}]^{n}
\leq C(1/\sigma)^{2\tau}\varepsilon_{0} .
$$

To prove the $(n+1)$-th step, it suffices to verify
the conditions in Lemma \ref{iteration lemma}. All together, we are led to showing the difference
 $\Omega h_{n}-\Omega h$
is small enough so that the non-degeneracy and invertibility
conditions hold. Furthermore, we also need to show
$DH_{n}$ and $(DH_{n})^{-1}$ are uniformly bounded along
all the iterations. These follow from \eqref{delta lambda} and \eqref{delta h} that
\begin{equation}\label{cauchy 2}
\begin{split}
 |\lambda_{n}-\lambda_{0}|, \|\Omega h_{n}-\Omega h_{0}\|_{\rho_{n}}
\leq &\sum_{j=1}^{n}\max\left\{|\Delta\lambda_{j-1}|, \|\Omega\Delta h_{j-1}\|_{\rho_{j}}\right\}\\
\leq &\sum_{j=1}^{\infty}C\dfrac{\varepsilon_{j-1}}{\sigma_{j}^{\tau}}
\leq
\sum_{j=1}^{\infty} \tilde\varepsilon_{j},
\end{split}
\end{equation}
provided that $DH_{0}$ is invertible and satisfies the non-degeneracy condition \eqref{KAM non-degeneracy}
and $\|e\|_{\rho}$ is sufficiently small.

Since $\rho_{n}$ decreases to $\rho/2$,
for the convergence of $h_{n}$ and $\lambda_{n}$
, it is sufficient
to apply the same estimates in \eqref{cauchy 2}
to show that $\{\Omega h_{n}\}_{n=0}^{\infty}$ is Cauchy on the uniform analyticity domain
 $\{x\in\TT^{d}_{\CC}=\CC^{d}/\ZZ^{d}: |\textrm{Im} x|\leq \rho/2 \}$, which is an immediate result of the convergence of
$\sum_{j=1}^{\infty}\tilde{\varepsilon}_{j}$.

This completes the proof of the convergence of the Newton iteration. \qed

\section*{Acknowledgement}
This work was initiated when X. H. was a visiting graduate student at Georgia Inst. of Technology. 
X. H. is grateful for the hospitality sponsored by CSC.


\appendix

\section{Phase locking phenomena for general
  torus maps}\label{general torus map}
In the Appendix, we summarize results for families
of maps close to rotations without the assumption of
preservation of a foliation. We hope that this can help
understand the effect of  the foliation preserving assumption
and see the effect in the applications where it appears.

There exist a lot of  results on the dynamical properties of torus maps both in
 mathematical and  physical literature.
 The systems with more than two frequencies have much more complicated dynamics than
the single frequency case. For instance, diffeomorphism on the torus can have
multiple attractors as well as chaotic trajectories.

In a generic family of torus maps close to
a family of rotations, the classical KAM theory asserts that the set of the parameter value for
which the diffeomorphism can not be smoothly conjugated to irrational rotations
occupies small Lebesgue measure in the parameter space.

There are also resonant regions in parameter space in which periodic orbits with given
rotation vector exist. Roughly speaking, resonance is a generalization of phase locking encountered
in circle maps with rational rotation number when the circle map has
an attracting or repulsive periodic point \cite{Arn88}. 
The resonant regions and phase locking phenomenon of dynamical systems on the torus have been studied
in \cite{KMG89, BGKM89, Gal89, Gal94, PdlLV03}.

In this Appendix, we  consider the phase locking phenomenon for a family of torus maps
close to a family of rotations on the torus. More precisely, we study
the following  families of torus map
\begin{equation*}
 F_{\Omega,\varepsilon}(x)=x+\Omega+\varepsilon f_{\varepsilon}(x),
\end{equation*}
where $x\in\torus^{d}$ and $(\Omega,\varepsilon)$ is
regarded as the parameter. We will study both the analytic case and finitely differentiable case
for  the torus map.

In contrast with \cite{Gal89,Gal94}, we
make further assumptions on the resonant frequency instead of restricting on the Mathieu type perturbations,
which enables us to obtain the convergence of the formal Fourier series.
Using averaging method and perturbation theory, we obtain a resonant normal form for
the torus maps and resonant regions
in the parameter space.

\subsection{Resonant normal form}\label{Normal Form general}
Consider
the following torus map
\begin{equation}\label{original system general}
 F_{\Omega, \varepsilon}(x)=x+\Omega+\varepsilon f_{\varepsilon}(x),
\end{equation}
where $x\in\torus^{d}$ and $f_{\varepsilon}$ is an analytic function defined on the torus for every $\varepsilon$.
We also assume some regularity of $f$ with respect to $\varepsilon$ such that the power series expansion in
$\varepsilon$ is valid.

We regard $\RR^{d}\times\RR^{+}$ of points $(\Omega,\varepsilon)$ as the parameter space.
Typically, when $\varepsilon=0$ and $\Omega$ is resonant,
the torus decomposes into a family of lower dimension torus. Each orbit is dense on such a lower dimension  torus, but not
in $\TT^{d}$.
We are interested in the dynamical behavior of $F_{\Omega, \varepsilon}$ for those parameter $\Omega$ close to a resonant frequency
and small $\varepsilon$.

In this section, we develop a resonant normal form for $F_{\Omega, \varepsilon}$ when $\Omega$ is close to the resonant frequency $\Omega_{0}$
whose resonance module is   $(d-r)$-dimension. From Subsection \ref{resonant intrinsic} we know that there exist a matrix $\mathfrak{A}\in \textrm{SL}(d, \ZZ)$,
a non-resonant intrinsic frequency
$\omega\in\RR^{r}$ and an integer vector $L\in\ZZ^{d}$ such that
\begin{equation}\label{intrinsic frequency general}
 \mathfrak{A}\Omega_{0}=
\begin{pmatrix}
 \omega\\
0
\end{pmatrix}
+L\ .
\end{equation}

Repeating the calculations in Subsection 2.1.1, we are able to obtain
the resonant normal form of $F_{\Omega,\varepsilon}$.
The main difference with the foliation preserving torus map is that the coordinate transformation $H_{\varepsilon}$ 
takes the form of $H_{\varepsilon}=Id+ \varepsilon h_{\varepsilon}$.

\begin{theorem}\label{lemma normal form general}
Given $p, q\in\NN$ with $p>q\geq (n+m-1)(d+\tau)$.
Assume $\Omega_{0}$ is resonant and the intrinsic frequency $\omega$
satisfies Diophantine condition \eqref{Diophantine condition}.
If $F_{\Omega, \varepsilon}(x)$ is real analytic $($or $\mathscr{C}^{p}$$)$ in $x$ and $\mathscr{C}^{n+m}$ in $\varepsilon$,
 then for sufficiently small
$\varepsilon$, the  torus
map $F_{\Omega,\varepsilon}$ is  conjugated by a real analytic $($or $\mathscr{C}^{p-q}$, respectively$)$ invertible
 function $$H_{\varepsilon}^{n+m-1}
 = Id+ \varepsilon h^{0}+ \varepsilon^{2} h^1+\cdots+\varepsilon^{n+m-1} h^{n+m-2} 
 $$ and a linear transformation
$\mathfrak{A}\in\textrm{SL}(d,\ZZ)$
to the following resonant normal form $\widetilde{F}_{\Omega, \varepsilon}: \TT^{r}\times\TT^{d-r}
 \rightarrow \TT^{r}\times\TT^{d-r} $ defined by
\begin{equation}\label{normal form equ general}
\widetilde{F}_{\Omega,\varepsilon}
\begin{pmatrix}x_{1}\\x_{2}\end{pmatrix}=
\begin{pmatrix}
x_{1}+\omega+\varepsilon^{n}\beta(x_{2},\varepsilon)+
\varepsilon^{n+m}r_{1}(x_{1},x_{2}; \varepsilon)
+\widetilde{\Delta}_{1}(\Omega)
\\
 x_{2}
+\varepsilon^{n}\eta(x_{2},\varepsilon)
+\varepsilon^{n+m}r_{2}(x_{1},x_{2}; \varepsilon)
+\widetilde{\Delta}_{2}(\Omega)
\end{pmatrix},
\end{equation}
where $n,m, h^{j}$ and  $\beta,\eta,r_{1},r_{2},\widetilde{\Delta}_{1}(\Omega),\widetilde{\Delta}_{2}(\Omega)$
can be computed explicitly from the averaging procedure.

Furthermore,  $r_{1}$ and $r_{2}$ are real analytic $($or $\mathscr{C}^{p-q}$, respectively$)$
periodic functions defined on $\TT^{d}$,
$\beta$ and $\eta$ are real analytic periodic $($or $\mathscr{C}^{p-q}$, respectively$)$ functions on $\TT^{d-r}$,
and $\widetilde{\Delta}_{j}(\Omega)=O(|\Omega-\Omega_{0}|), j=1,2$.
\end{theorem}

\subsection{Persistence of the invariant surface}\label{Persistence} By an abuse of notation, we still
denote by $F_{\Omega, \varepsilon}$  the resonant normal form developed in Theorem \ref{lemma normal form general}.
In this section, we assume that there exists $x_{2}^*\in\TT^{d-r}$  such that
\begin{equation}\tag{H1}\label{zero point assumption general}
\eta(x_{2}^*,0)=0\  \textrm{and}\  D_{1}\eta(x_{2}^{*},0)\  \textrm{is non-singular}.
\end{equation}
By the implicit function theorem,
we obtain $x_{2}^{*}(\varepsilon)$ such that $\eta(x_{2}^{*}(\varepsilon),\varepsilon)\equiv 0$
 and $x^{*}_{2}(0)=x_{2}^*$.
Then it follows that $\Gamma=\{(x_{1},x_{2}^{*}(\varepsilon))\ |\ x_{1}\in\TT^{r}\}$ is an invariant surface
for the torus map
$(x_{1},x_{2})\mapsto(x_{1}+\omega+\varepsilon^{n}\beta(x_{2},\varepsilon),x_{2}
+\varepsilon^{n}\eta(x_{2},\varepsilon))$
on which the motion is a rotation.
Without loss of generality, we assume $x_{2}^{*}(\varepsilon)=0$.

In what follows, we show the existence of invariant surface for the resonant normal
form developed in Theorem \ref{lemma normal form general}, which is close to $\Gamma$.

\begin{theorem}\label{phase locking theorem gen}
 Let $F_{\Omega,\varepsilon}$ be a family of real analytic ($\mathscr{C}^{\infty}$
or finitely differentiable) torus maps given by \eqref{original system general}.
Assume $\Omega_{0}$ is resonant with intrinsic frequency $\omega$ satisfying
Diophantine condition \eqref{Diophantine condition}.
 Then if
the resonant normal form \eqref{normal form equ general} satisfies non-degeneracy condition \eqref{zero point assumption general},
 there exists a $r$-dimensional,  finitely differentiable invariant torus
of $F_{\Omega, \varepsilon}$ for  those parameters $(\Omega,\varepsilon)$ in some neighborhood
of $(\Omega_{0},0)$.

In the finitely differentiable case (but with enough regularity), the regularities of the low dimensional invariant tori
are less than those of $F_{\Omega,\varepsilon}$.
\end{theorem}

In the perturbative setting, we see that the invariant manifold
of $F_{\Omega,\varepsilon}$ can be represented by a graph
of $w: \TT^{r}\rightarrow\RR^{d-r}$, whose invariance determines a functional $\mathscr{F}$ defined
on some function space. Then by the contraction mapping arguments, we prove the existence
of a fixed point for $\mathscr{F}$, which corresponds to the desired invariant manifold of $F_{\Omega,\varepsilon}.$
We delay the proof of Theorem \ref{phase locking theorem gen} to the end of this section.

\subsubsection{Formulation of the functional $\mathscr{F}$}
For the resonant normal form for $F_{\Omega,\varepsilon}$, we further assume
the hyperbolicity that

\begin{equation}\label{hyperbolic without pert.}
 \textrm{Spec}(I+\varepsilon^{n}D_{1}\eta(0,\varepsilon))\cap\{z\in\ZZ\ :\ |z|=1\}=\emptyset
\end{equation}
holds for small $\varepsilon>0$.

By the hyperbolicity assumption \eqref{hyperbolic without pert.}, we
decompose the space $\RR^{d-r}$ into $\RR^{d-r}=E^{s}_{\varepsilon}\oplus E^{u}_{\varepsilon}$
 such that
\begin{equation*}
 I+\varepsilon^{n} D_1\eta(0,\varepsilon)=
\begin{pmatrix}
 \Lambda^{s}_{\varepsilon} & 0\\
0& \Lambda^{u}_{\varepsilon}
\end{pmatrix},
\end{equation*}
where the spectral radius of $\Lambda^{s}_{\varepsilon}$ is less than one and $\Lambda^{u}_{\varepsilon}$
is invertible with spectral radius greater than one.

Assume that there exist positive constants $C_{s}$, $C_{u}$
and adapted norms for $E^{s}_{\varepsilon}$ and $E^{u}_{\varepsilon}$
such that
\begin{equation}\tag{H2}\label{adapted norms}
 \|\Lambda^{s}_{\varepsilon}\|\leq 1-C_{s}\varepsilon^{n}\quad
\textrm{and}\quad \|(\Lambda^{u}_{\varepsilon})^{-1}\|\leq 1-C_{u}\varepsilon^{n}.
\end{equation}

\begin{remark}\label{spec radius foliation}
We know that,
for the $\Omega$-foliation preserving torus maps $($see \eqref{original system} and Theorem
\ref{lemma normal form}$)$,  $\eta(x_{2},\varepsilon)$ is a scalar
function. As a result, assumption \eqref{zero point assumption general} implies \eqref{adapted norms}
for  small $\varepsilon$. More precisely, we can always choose constants $C_{s}$
and $C_{u}$ such that
$$
1+\varepsilon^{n}D_{1}\eta(0,\varepsilon)\leq 1-C_{s}\varepsilon^{n}\quad
\textrm{if}\quad D_{1}\eta(0,0)<0,
$$
and
$$
(1+\varepsilon^{n}D_{1}\eta(0,\varepsilon))^{-1}\leq 1-C_{u}\varepsilon^{n}\quad
\textrm{if}\quad D_{1}\eta(0,0)>0,
$$
for sufficiently small $\varepsilon$.
In this case, it suffices to consider two different situations rather than
to decompose the space.
\end{remark}

In the perturbative setting, we assume that the invariant manifold
of $F_{\Omega,\varepsilon}$, which is close to $\Gamma$, can be represented by a graph
of $w: \TT^{r}\rightarrow\RR^{d-r}$.
Then by the invariance of the graph, we are led to solving the following functional equation
\begin{equation}\label{graph invariance general}
\begin{split}
 & [I+\varepsilon^{n}D_{1}\eta(0,\varepsilon)]w(x_{1})+\varepsilon^{n+m}r_{2}(x_{1},w(x_{1}))+\varepsilon^{n}R(w(x_{1}),\varepsilon)
 +\widetilde{\Delta}_{2}(\Omega)\\
 =& w\left(x_{1}+\omega+\varepsilon^{n}\beta(w(x_{1}),\varepsilon)+\varepsilon^{n+m}r_{1}(x_{1},w(x_{1}))
+\widetilde{\Delta}_{1}(\Omega)\right)
 \end{split}
\end{equation}
for every $x_{1}\in\TT^{r}$,
where
\begin{equation*}
 R(x_{2},\varepsilon)=\eta(x_{2},\varepsilon)-D_{1}\eta(0,\varepsilon)x_{2}
=\int_{0}^{1}\int_{0}^{1}sD_{11}\eta(stx_{2},\varepsilon)\mathrm{d}s\mathrm{d}t\cdot x_{2}^{\otimes 2}.
\end{equation*}
Projecting \eqref{graph invariance general} onto $E^{s}_{\varepsilon}$ and $E^{u}_{\varepsilon}$ respectively,
we obtain that \eqref{graph invariance general} is equivalent to
\begin{equation}\label{stable subspace}
 \Lambda^{s}_{\varepsilon}w^{s}+\varepsilon^{n+m} r_{2}^{s}
 +\varepsilon^{n}
R^{s}
+\widetilde{\Delta}_{2}^{s}(\Omega)
=w^{s}\circ \mathcal{P}
\end{equation}
and
\begin{equation}\label{unstable subspace}
 \Lambda^{u}_{\varepsilon}w^{u}+\varepsilon^{n+m} r_{2}^{u}
 +\varepsilon^{n}
R^{u}
+\widetilde{\Delta}_{2}^{u}(\Omega)
=w^{u}\circ\mathcal{P},
\end{equation}
where
$$w=(w^{s},w^{u}),\ r_{2}=(r_{2}^{s},\ r_{2}^{u}),\ R=(R^{s},R^{u}),\
\widetilde{\Delta}_{2}(\Omega)=(\widetilde{\Delta}_{2}^{s}(\Omega),\widetilde{\Delta}_{2}^{u}(\Omega)),$$
and
\begin{align*}
& r_{2}^{s}
=r_{2}^{s} (x_{1},w^{s}(x_{1}),w^{u}(x_{1})),\quad
r_{2}^{u}= r_{2}^{u}
 (x_{1},w^{s}(x_{1}),w^{u}(x_{1})),
\\
& R^{s}
=R^{s} (w^{s}(x_{1}),w^{u}(x_{1}),\varepsilon),\quad
R^{u}
=R^{u}(w^{s}(x_{1}),w^{u}(x_{1}),\varepsilon),\\
&\mathcal{P}(x_{1})=x_{1}+\omega+\varepsilon^{n}\beta(w(x_{1}),\varepsilon)+\varepsilon^{n+m}r_{1}(x_{1},w(x_{1}))
+\widetilde{\Delta}_{1}(\Omega).
\end{align*}

Now we are looking for a functional $\mathscr{F}$ such that its fixed point is a solution
 to \eqref{stable subspace} and \eqref{unstable subspace}. Before this,
for \eqref{stable subspace}, we are led to solving $x_{1}$ from
\begin{equation}\label{priori equ general}
 y_{1}=\mathcal{P}(x_{1})=x_{1}+\omega+\varepsilon^{n}\beta(w(x_{1}),\varepsilon)+\varepsilon^{n+m}r_{1}(x_{1},w(x_{1}))
+\widetilde{\Delta}_{1}(\Omega).
\end{equation}
Similar arguments also appear in \cite{Lan73}.

Rewrite \eqref{priori equ general} as
\begin{equation*}
 \begin{split}
  x_{1}=&y_{1}-\omega-\varepsilon^{n}\beta(w^{s}(x_{1}),w^{u}(x_{1}),\varepsilon)-\varepsilon^{n+m}
r_{1}(x_{1},w^{s}(x_{1}),w^{u}(x_{1}))-\widetilde{\Delta}_{1}(\Omega)\\
\equiv&\mathscr{G}(x_{1};y_{1},w^{s},w^{u}),
 \end{split}
\end{equation*}
and consider the fixed point problem of operator $\mathscr{G}$ with parameter $(y_{1},w)$
belonging to $\TT^{r}\times \mathcal{U}$, where
$$\mathcal{U}= \{w^{s}\in\mathscr{C}^{1}(\TT^{r},E^{s}):\ \|w^{s}\|_{1}\leq 1\}
\times \{w^{u}\in\mathscr{C}^{1}(\TT^{r},E^{u}):\ \|w^{u}\|_{1}\leq 1\}.$$
Noticing that
\begin{align*}
 |\mathscr{G}(x_{1}^{(1)};y_{1},w^{s},w^{u})&-\mathscr{G}(x_{1}^{(2)};y_{1},w^{s},w^{u})|\\
 &\leq
\varepsilon^{n}\left(\|D\beta\|_{0}+2\varepsilon^{m}\|Dr_{1}\|_{0}+C\frac{|\Delta\Omega|}{\varepsilon^{n}}\right)\cdot |x_{1}^{(1)}-x_{1}^{(2)}|,
\end{align*}
we have for
\begin{equation}\label{Delta Omega small}
 |\Delta\Omega|\leq C\varepsilon^{n+1}
\end{equation}
and sufficiently small $\varepsilon$,
the operator is a uniform contraction mapping on $\TT^{r}\times\mathcal{U}$.

By the uniform contraction mapping theorem (see \cite{chi06}), we obtain the solution
of \eqref{priori equ general} denoting by $x_{1}=v(y_{1},w^{s},w^{u})$.
Furthermore,
since $\mathscr{G}:  \TT^{r}\times \mathcal{U}\rightarrow \TT^{r}$ is $\mathscr{C}^{1}$,
we have the fixed point function $v: \TT^{r}\times\mathcal{U}\rightarrow \TT^{r}$ is also $\mathscr{C}^1$.

Then we rewrite \eqref{stable subspace} as
\begin{equation*}
 \begin{split}
  w^{s}(y_{1})=&\Lambda_{\varepsilon}^{s}w^{s}(v(y_{1},w^{s},w^{u}))
  +\varepsilon^{n+m}r_{2}^{s}(v,w^{s}(v),w^{u}(v)) \\
  & +
\varepsilon^{n}R^{s}(w^{s}(v),w^{u}(v),\varepsilon)+\widetilde{\Delta}^{s}_{2}(\Omega)\\
\equiv & \mathscr{F}^{s}[w^{s},w^{u}](y_{1}).
 \end{split}
\end{equation*}
For \eqref{unstable subspace}, it would be much easier to find a functional $\mathscr{F}^{u}$ whose
fixed point together with $\mathscr{F}^{s}$ would be a solution of \eqref{stable subspace} and
\eqref{unstable subspace}. Actually, we have
\begin{equation*}
\begin{split}
 w^{u}(x_{1})=&(\Lambda^{u}_{\varepsilon})^{-1}w^{u}(\mc{P}(x_{1}))
-\varepsilon^{n+m}(\Lambda_{\varepsilon}^{u})^{-1}r^{u}_{2}(x_{1}, w^{s}(x_{1}),w^{u}(x_{1}))\\
&-\varepsilon^{n}(\Lambda^{u})^{-1}R^{u}(w^{s}(x_{1}),w^{u}(x_{1}),\varepsilon)
-(\Lambda^{u})^{-1}\widetilde{\Delta}_{2}^{u}(\Omega)\\
\equiv& \mathscr{F}^{u}[w^{s},w^{u}](x_{1}).
\end{split}
\end{equation*}

All together, we define
\begin{equation}
 \mathscr{F}[w^{s},w^{u}]=
(\mathscr{F}^{s}[w^{s},w^{u}], \mathscr{F}^{u}[w^{s},w^{u}]\ ).
\end{equation}
from $\mathscr{C}^{L}(\TT^{r},E_{\varepsilon}^{s})\times \mathscr{C}^{L}(\TT^{r},E_{\varepsilon}^{u})$ to itself.

It is readily seen that, given an integer $L>0$ and a real number $0<l\leq n$, there
exists a $\varepsilon_{0}>0$ such that for all $0<\varepsilon<\varepsilon_{0}$, we have
\begin{equation}\label{small condition s}
\|\Lambda^{s}_{\varepsilon}\|(1-\varepsilon^{n+l}\|D\beta\|_{0})^{-j}<1
\end{equation}
and
\begin{equation}\label{small condition u}
 \|(\Lambda^{u}_{\varepsilon})^{-1}\|(1+\varepsilon^{n+l}\|D\beta\|_{0})^{j}<1
\end{equation}
for any $0\leq j\leq L$.

\begin{remark}
The adapted norms of $\Lambda^{s}$ and $(\Lambda^{u})^{-1}$ depend also
on the perturbation. We will see more clearly about the smallness of $\varepsilon$ in the
one dimensional resonance case. For instance, in one dimensional resonance case and
$D_{1}\eta(0,0)>0$, condition \eqref{small condition u} reduces to
\begin{equation*}
 \frac{1}{1+\varepsilon^{n} D_{1}\eta(0,0)}(1+\varepsilon^{n+l}\|D\beta\|_{0})^{L}<1.
\end{equation*}
Then, for sufficiently large $L$, we see that $\varepsilon\leq C\sqrt[l]{\frac{1}{L}}$.
\end{remark}

\subsubsection{Invariant closed subset of the functional}
Now we are trying to find a closed subset $\mathscr{L}\subset \mathscr{C}^{L}(\TT^{r},E^{s})
\times \mathscr{C}^{L}(\TT^{r},E^{u})$ such that $\mathscr{F}(\mathscr{L})\subseteq\mathscr{L}$.
In this paper, we would like to find the subset $\mathscr{L}$ with different bounds for
different orders of derivatives. More precisely, we assume
\begin{equation}\label{form of L}
\begin{split}
 \mathscr{L}=\Big\{(w^{s},w^{u})&\in\mathscr{C}^{L}(\TT^{r},E^{s})\times
 \mathscr{C}^{L}(\TT^{r},E^{u}):\\
 &
\|D^{j}w^{s}\|_{0}\leq \delta_{j}, \|D^{j}w^{u}\|_{0}\leq \delta_{j}~
 \textrm{for all}\ 0\leq j\leq L\Big\}.
\end{split}
\end{equation}
In what follows, we show how to find $\delta_{j}$'s such that $\mathscr{F}(\mathscr{L})\subseteq
\mathscr{L}$.

For $\|w^{s}\|_{0}\leq \delta_{0}<1$ and  $\|w^{u}\|_{0}\leq \delta_{0}<1$, we have
\begin{equation*}
\begin{split}
 &\|\mathscr{F}^{s}[w^{s},w^{u}]\|_{0}\\
 \leq& \|\Lambda_{\varepsilon}^{s}\|\delta_{0}+\varepsilon^{n+m}\|r_{2}^{s}\|_{0}
 +\varepsilon^{n}\|\eta(\cdot,0)\|_{2}\delta^{2}_{0}+\varepsilon^{n+1}\|D_{2}\eta(\cdot, 0)\|_{2}
 +C|\Delta\Omega|\\
 \leq &
 (1-C_{s}\varepsilon^{n})\delta_{0}+\varepsilon^{n+1}\|D_{2}\eta(\cdot, 0)\|_{2}+\varepsilon^{n+m}\|r_{2}^{s}\|_{0}
 +\varepsilon^{n}\|\eta(\cdot, 0)\|_{2}\delta^{2}_{0}+C|\Delta\Omega|.
\end{split}
\end{equation*}
To ensure $\|\mathscr{F}^{s}[w^{s}, w^{u}]\|_{0}\leq \delta_{0}$, it suffices
\begin{equation*}
\|\eta\|_{2}\delta_{0}^{2}-C_{s}\delta_{0}+\varepsilon \Bigg(\|\eta\|_{1}+\varepsilon^{m}\|\eta\|_{2}+C
|\Delta\Omega|\varepsilon^{-n}\Bigg)\leq 0,
\end{equation*}
which implies
\begin{equation}\label{range delta 0}
O(\varepsilon)=C_{s}\frac{1-\sqrt{1-4\varepsilon\|\eta\|_{2}q}}{2\|\eta\|_{2}}\leq \delta_{0}\leq C_{s} \frac{1+\sqrt{1-4\varepsilon\|\eta\|_{2}q}}{2\|\eta\|_{2}}=O(1)
\end{equation}
with $q=\|\eta\|_{1}+\varepsilon^{m}\|\eta\|_{2}+C|\Delta\Omega|\varepsilon^{-n}$.
Then we see that, when $0<\alpha_{0}<1$, $\delta_{0}=\varepsilon^{\alpha_{0}}$
satisfies \eqref{range delta 0} for sufficiently small $\varepsilon$.

For $\mathscr{F}^{u}[w^{s},w^{u}]$, we also have
\begin{equation*}
\begin{split}
 &\|\mathscr{F}^{u}[w^{s},w^{u}]\|_{0}\\
 \leq& \|(\Lambda^{u})^{-1}\|\delta_{0}+\varepsilon^{n+m}\|r_{2}^{s}\|_{0}
 +\varepsilon^{n}\|\eta\|_{2}\delta^{2}_{0}+\varepsilon^{n+1}\|\eta\|_{1}
 +C|\Delta\Omega|\\
 \leq &
 (1-C_{u}\varepsilon^{n})\delta_{0}+\varepsilon^{n+1}\|\eta\|_{1}+\varepsilon^{n+m}\|r_{2}^{s}\|_{0}
 +\varepsilon^{n}\|\eta\|_{2}\delta^{2}_{0}+C|\Delta\Omega|.
\end{split}
\end{equation*}
Similarly, it suffices to choose $\delta_{0}=\varepsilon^{\alpha_{0}}$ with $0<\alpha_{0}<1$ such that
$$\|\mathscr{F}^{u}[w^{s},w^{u}]\|_{0}\leq \delta_{0}.$$

In what follows, we also need to estimate the derivatives of $\mathscr{F}[w^{s},w^{u}]$ up to the
order of $L$. We will show that the bound for the first derivative of $\mathcal{P}$ plays an important role
in the analysis of other higher order derivatives.

The derivative of $\mathscr{F}^{u}[w]$ is
\begin{equation*}
\begin{split}
 &D(\mathscr{F}^{u}[w^{s},w^{u}])(x_{1})\\
=&
(\Lambda^{u}_{\varepsilon})^{-1}Dw^{u}(\mathcal{P}(x_{1}))\cdot \Bigg\{
I+\varepsilon^{n}D_{1}\beta(\cdot)\cdot \Big[Dw^{s}(x_{1}),Dw^{u}(x_{1})\Big]\\
&+\varepsilon^{n+m}D_{1}r_{1}(\cdot)+
\varepsilon^{n+m}D_{2}r_{1}(\cdot)\Big[Dw^{s}(x_{1}), Dw^{u}(x_{1})\Big]+D_{x_{1}}[\widetilde{\Delta}_{1}(\Omega)]\Bigg\}\\
&-\varepsilon^{n+m}(\Lambda^{u}_{\varepsilon})^{-1} D_{1}r_{2}^{u}(\cdot)-\varepsilon^{n+m}
(\Lambda^{u}_{\varepsilon})^{-1} D_{2}r_{2}^{u}(\cdot)\Big[Dw^{s}(x_{1}), Dw^{u}(x_{1})\Big]\\
&-\varepsilon^{n}(\Lambda_{\varepsilon}^{u})^{-1}D_{1}R^{u}(\cdot)\Big[Dw^{s}(x_{1}),Dw^{u}(x_{1})\Big]
-(\Lambda_{\varepsilon}^{u})^{-1} D_{x_{1}}[\widetilde{\Delta}_{2}^{u}(\Omega)],
\end{split}
\end{equation*}
where
$(\cdot)$ denotes the argument omitted in a function.

Similarly, the derivative of $\mathscr{F}^{s}[w^{s},w^{u}]$ reads
 \begin{equation*}
\begin{split}
  & D(\mathscr{F}^{s}[w^{s},w^{u}])(y_{1})\\
  =&\Lambda_{\varepsilon}^{s} Dw^{s}(v)D_{1}v(\cdot)+
\varepsilon^{n+m}\Bigg\{D_{1}r_{2}^{s}(\cdot)+D_{2}r_{2}^{s}\Big[Dw^{s}(v),Dw^{u}(v)\Big]
\Bigg\}\cdot D_{1}v(\cdot)\\
&+\varepsilon^{n}D_{1}R^{s}(\cdot)\Big[Dw^{s}(v),Dw^{u}(v)\Big]\cdot D_{1}v(\cdot)
+D_{y_{1}}[\widetilde{\Delta}_{2}^{s}(\Omega)].
\end{split}
\end{equation*}
However, the derivative of $\mathscr{F}^{s}$ is more involved. We also need to give the derivative
of $v(y_{1},w^{s},w^{u})$ with respect to $y_{1}$. From \eqref{priori equ general} we obtain
\begin{align*}
  D_{1}v(y_{1}, w^{s},w^{u})=&\Bigg\{I+\varepsilon^{n}D\beta(\cdot)\Big[Dw^{s}(v),Dw^{u}(v)\Big]
+\varepsilon^{n+m}D_{1}r_{1}(\cdot)\\
&+\varepsilon^{n+m}D_{2}r_{1}(\cdot)\Big[Dw^{s}(v),Dw^{u}(v)\Big]
+D_{x_{1}}[\widetilde{\Delta}_{1}(\Omega)]
\Bigg\}^{-1}.
\end{align*}

For $\|Dw^{s}\|_{0}\leq \delta_{1}$ and $\|Dw^{u}\|_{0}\leq \delta_{1}$, we have
\begin{align*}
 & \|D\mathscr{F}^{u}[w^{s},w^{u}]\|_{0}\\
\leq &(1-C_{u}\varepsilon^{n})
\Bigg\{
  \delta_{1}\Big[1+\varepsilon^{n}\|D\beta\|_{0}\delta_{1}+\varepsilon^{n+m}\|r_{1}\|_{1}(1+\delta_{1})
+C|\Delta\Omega|\Big]\\
 & +\varepsilon^{n+m}\|r_{2}\|_{1}(1+\delta_{1})
  +\varepsilon^{n+1}\|\eta\|_{1}\delta_{1}+\varepsilon^{n}\|\eta\|_{2}\delta_{0}\delta_{1}
+C|\Delta\Omega|\Bigg\}.\\
\leq &
(1-C_{u}\varepsilon^{n})(1+\varepsilon^{n}\|D\beta\|_{0}\delta_{1})\delta_{1}
+\varepsilon^{n+\alpha_{0}}\|\eta\|_{2}+\varepsilon^{n+1}\left(\|\eta\|_{1}+C |\Delta\Omega|\varepsilon^{-(n+1)}\right)
\\
&+2\varepsilon^{n+m}(\|r_{1}\|_{1}+\|r_{2}\|_{1}).
\end{align*}
Then it suffices to choose $\delta_{1}=\varepsilon^{\alpha_{1}}$ with $l<\alpha_{1}<n+\alpha_{0}$
such that $\|D\mathscr{F}^{u}[w^{s},w^{u}]\|_{0}\leq \delta_{1}$.

To estimate the derivative of $\mathscr{F}^{s}[w^{s},w^{u}]$, we see that
\begin{equation*}
  \|D_{1}v(\cdot,w^{s},w^{u})\|_{0}
\leq  (1-\varepsilon^{n}\|D\beta\|_{0}\delta_{1})^{-1}+2\|r_{1}\|_{1}\varepsilon^{n+m}+C|\Delta\Omega|
\end{equation*}
for sufficiently small $\varepsilon$.

Consequently, we  have
\begin{equation*}
 \begin{split}
  \|D\mathscr{F}^{s}[w^{s},w^{u}]&\|_{0}
\leq (1-C_{s}\varepsilon^{n})(1-\varepsilon^{n}\|D\beta\|_{0}
\delta_{1})^{-1}\delta_{1}+4\varepsilon^{n+\alpha_{0}}\|\eta(\cdot,0)\|_{2}\\
+&4\varepsilon^{n+1}\left(\|D_{2}\eta(\cdot,0)\|_{1}+C|\Delta\Omega|\varepsilon^{-(n+1)}\right)
+4\varepsilon^{n+m}(\|r_{1}\|_{1}+\|r_{2}\|_{1}),
 \end{split}
\end{equation*}
which implies that the choice of $\delta_{1}=\varepsilon^{\alpha_{1}}$ with $l<\alpha_{1}<n+\alpha_{0}$
is sufficient to keep $\|D\mathscr{F}^{s}[w^{s},w^{u}]\|_{0}\leq \delta_{1}$.

For higher derivatives, the Faa-di-Bruno's formula (see \cite[p.3]{AR67}) gives
\begin{equation}\label{Faa-di-Bruno u}
 \begin{split}
  &~ D^{i}\mathscr{F}^{u}[w^{s},w^{u}](x_{1})\\
=&  (\Lambda^{u})^{-1}D^{i}
w^{u}(\mathcal{P}(x_{1}))\cdot [
D\mathcal{P}]^{\otimes i}\\
&+(\Lambda^{u})^{-1}\sum_{1\leq q<i}\sum_{k_{1}+\cdots+k_{q}=i}
\sigma_{i}D^{q}w^{u}\circ \mathcal{P}\ \Big[D^{k_1}\mathcal{P},\cdots,D^{k_q}\mathcal{P}\Big]\\
&-\varepsilon^{n+m}(\Lambda^{u})^{-1}D^{i}(r_{2}(x_{1},w(x_{1})))
-\varepsilon^{n}(\Lambda^{u})^{-1}\\
&~\times
\sum_{1< q\leq i}\sum_{k_{1}+\cdots+k_{q}=i}
\sigma_{i}D^{q}\eta^{u}\circ w\ \Big[(D^{k_1}w^{s},D^{k_1}w^{u}),\cdots, (D^{k_q}w^{s},D^{k_q}w^{u})\Big],
 \end{split}
\end{equation}
where
\begin{align*}
D\mathcal{P}=&I+\varepsilon^{n}D_{1}\beta(\cdot)\cdot [Dw^{s}(x_{1}),Dw^{u}(x_{1})]
+\varepsilon^{n+m}D_{1}r_{1}(\cdot)\\
&+\varepsilon^{n+m}D_{2}r_{1}(\cdot)[Dw^{s}(x_{1}), Dw^{u}(x_{1})]
+D_{x_1}[\widetilde{\Delta}_{1}(\Omega)],
\end{align*}
$k_{1},\cdots,k_{q}$ are positive integers
and $\sigma_{i+1}=\sigma_{i}(k_{1},\cdots,k_{q})$ is of integer value
and can be calculated explicitly.

Noticing that when $1\leq q<i$, there always exists a $k_{j}$ with $1\leq j\leq q$ such that $k_{j}\geq 2$.
Thus $D^{k_{j}}\mathcal{P}(x_{1})$ will be of order $\varepsilon^{n}$.

By induction, we assume that for $2\leq j\leq i-1$, $\delta_{j}=\varepsilon^{\alpha_{j}}$ with
$0<\alpha_{j}\leq n$. Then, we see
from \eqref{Faa-di-Bruno u} that
\begin{equation}\label{high derivative estimate}
\begin{split}
 \|D^{i}\mathscr{F}^{u}[w^{s},w^{u}]\|_{0}
\leq &(1-C_{u}\varepsilon^{n})(1+\varepsilon^{n+l}
\|D\beta\|_{0})^{L}\delta_{i}\\
&+\varepsilon^{n}\mathfrak{R}^{u}_{i}(\delta_{0},\cdots, \delta_{i-1};\varepsilon,
\eta,r_{1},r_{2},L)+C|\Delta\Omega|,
\end{split}
\end{equation}
where $\mathcal{R}^{u}_{i}$ is a polynomial of $\delta_{0},\cdots,\delta_{i-1}$ whose coefficients
depend on $\eta$, $r_{1}$, $r_{2}$, $L$ and $\varepsilon$.
By assumption \eqref{small condition u}, we have $\delta_{i}=\varepsilon^{\alpha_{i}}$ with
$0<\alpha_{i}\leq n$ is sufficient to ensure
$\|D^{i}\mathscr{F}^{u}[w^{s},w^{u}]\|_{0}\leq \delta_{i}$ when $\varepsilon$ is small.
This completes the induction argument.

For the high derivatives of $\mathscr{F}^{s}[w^{s},w^{u}]$, one also have the similar
estimate to \eqref{high derivative estimate} as follows
\begin{equation*}
\begin{split}
 \|D^{i}\mathscr{F}^{s}[w^{s},w^{u}]\|_{0}\leq& (1-C_{s}\varepsilon^{n})(1-\varepsilon^{n+l}
\|D\beta\|_{0})^{-L}\delta_{i}\\
&+
\varepsilon^{n}\mathfrak{R}^{s}_{i}(\delta_{0},\cdots, \delta_{i-1};\varepsilon,
\eta,r_{1},r_{2},L) +C|\Delta\Omega|,
\end{split}
\end{equation*}
which implies that the $\alpha_{j}$'s obtained keeps $\|D^{j}\mathscr{F}^{s}[w^{s},w^{u}]\|_{0}\leq \delta_{j}$.

We conclude the above arguments in the following lemma.
\begin{lemma}\label{invariant subset L}
Let $L\in\NN$ be fixed  and real number $l$ satisfy $0<l\leq n$.
 Let $F_{\Omega, \varepsilon}\in\mathscr{C}^{p}(\TT^{d}, \TT^{d})$ be the torus map in the resonant normal form \eqref{normal form equ}
with $L+2<p\in\NN$.
Assume \eqref{zero point assumption} and \eqref{hyperbolic without pert.} hold.
Then
for any given $\alpha_{j}$ satisfying $0<\alpha_{0}<1$, $l<\alpha_{1}<n+\alpha_{0}$
and $0<\alpha_{i}\leq n$ for $2\leq i\leq L$,
 there exists a $\varepsilon_{0}=\varepsilon_{0}(L,l,\alpha_{0},\cdots,\alpha_{L};\eta,r_{1},r_{2})$ such
that for any $0<\varepsilon<\varepsilon_{0}$ and $|\Omega-\Omega_{0}|\leq \varepsilon^{n+1}$, we have
$\mathscr{F}(\mathscr{L})\subseteq \mathscr{L}$, where $\mathscr{L}$ is given by \eqref{form of L}
with $\delta_{j}=\varepsilon^{\alpha_{j}}$.
\end{lemma}

\begin{remark}
 Since $L$ is given and fixed, one should keep in mind that
 there are only finitely many conditions on the smallness of $\varepsilon$
which appear in determining $\delta_{j}$ for $0\leq j\leq L$.
\end{remark}

\subsubsection{Contractility of the functional}
In this subsection, we show that $\mathscr{F}$ is a contraction on $\mathscr{L}$ in the $\mathscr{C}^{0}$-norm.

For $(w^{s}_{1},w^{u}_{1})$ and $(w^{s}_{2},w^{u}_{2})$ in $\mathscr{L}$, we have
\begin{align*}
 \|\mathscr{F}^{u}&[w^{s}_{2},w^{u}_{2}]-\mathscr{F}^{u}[w^{s}_{1},w^{u}_{1}]\|_{0}\\
\leq& \|(\Lambda^{u})^{-1}\|\cdot\Bigg\{
\|w^{u}_{2}-w^{u}_{1}\|_{0}+\|Dw^{u}_{1}\|_{0}\cdot\Big[
\varepsilon^{n}\|D\beta\|_{0}\cdot\|(w^{s}_{2}-w^{s}_{1},w^{u}_{2}-w^{u}_{1})\|_{0}\\
&+\varepsilon^{n+m}\|Dr_{1}\|_{0}\cdot\|(w^{s}_{2}-w^{s}_{1},w^{u}_{2}-w^{u}_{1})\|_{0}\Big]
+\varepsilon^{n+m}\|Dr_{2}\|_{0}\\
&\times\|(w^{s}_{2}-w^{s}_{1},w^{u}_{2}-w^{u}_{1})\|_{0}
+\varepsilon^{n}\|\eta\|_{2}\cdot (\|(w^{s}_{2},w^{u}_{2})\|_{0}+
\|(w^{s}_{1},w^{u}_{1})\|_{0})\\
&\times\|(w^{s}_{2}-w^{s}_{1},w^{u}_{2}-w^{u}_{1})\|_{0}+\varepsilon^{n}\|\eta\|_{3}\cdot
\|(w^{s}_1,w^{u}_1)\|_{0}^{2}\cdot\|(w^{s}_{2}-w^{s}_{1},w^{u}_{2}-w^{u}_{1})\|_{0}\\
&+\varepsilon^{n+1}\left(\|\eta\|_{1}+C |\Delta\Omega|\varepsilon^{-(n+1)}\right)\cdot
\|(w^{s}_{2}-w^{s}_{1},w^{u}_{2}-w^{u}_{1})\|_{0}
\Bigg\}\\
\leq &\Bigg\{(1-C_{u}\varepsilon^{n})(1+\varepsilon^{n+\alpha_{1}}\|D\beta\|_{0})+\varepsilon^{n+m+\alpha_{1}}
\|Dr_{1}\|_{0}+\varepsilon^{n+m}\|Dr_{2}\|_{0}\\
&+2\varepsilon^{n+\alpha_{0}}\|\eta\|_{2}+
\varepsilon^{n+2\alpha_{0}}+\varepsilon^{n+1}\left(\|\eta\|_{1}+C\frac{|\Delta\Omega|}{\varepsilon^{n+1}}\right)\Bigg\}
\cdot\|(w^{s}_{2}-w^{s}_{1},w^{u}_{2}-w^{u}_{1})\|_{0},
\end{align*}
which enables us to choose sufficiently small $\varepsilon$ such that $\mathscr{F}^{u}$
is a contraction in the $\mathscr{C}^{0}$-norm.

To show a similar estimate for $\mathscr{F}^{s}$, we give  the Lipschitz property of
$x_{1}=v(y_{1},w^{s},w^{u})$ with respect to $(w^{s},w^{u})$ on $\mathcal{U}$.
Obviously, from \eqref{priori equ general}, we see that
\begin{align*}
  -(v(y_{1},w_2)&-v(y_{1},w_{1}))\\=&\varepsilon^{n}\Big[\beta(w_{2}(v(y_{1},w_2)))-\beta(w_{1}(v(y_{1},w_2))\Big]\\
&+\varepsilon^{n}\Big[\beta(w_{1}(v(y_{1},w_2)))-\beta(w_{1}(v(y_{1},w_1))\Big]\\
&+\varepsilon^{n+m}\Big[r_{1}(v(y_1,w_2),w_2(v(y_1,w_2))-r_1(v(y_1,w_1),w_2(v(y_1,w_2))\Big]\\
&+\varepsilon^{n+m}\Big[r_{1}(v(y_1,w_1),w_2(v(y_1,w_2))-r_1(v(y_1,w_1),w_2(v(y_1,w_1))\Big]\\
&+\varepsilon^{n+m}\Big[r_{1}(v(y_1,w_1),w_2(v(y_1,w_1))-r_1(v(y_1,w_1),w_1(v(y_1,w_1))\Big]\\
&+\widetilde{\Delta}_{1}(\Omega; w_{2})-\widetilde{\Delta}_{1}(\Omega; w_{1}),
\end{align*}
where $w_{2}=(w_{2}^{s},w_{2}^{u})$ and $w_{1}=(w_{1}^{s},w_{1}^{u})$.

For $w^{1}, w^{2}\in\mathscr{L}$, one readily sees
\begin{equation*}
 \begin{split}
  |v(y_{1},w_2)-v(y_1,w_1)|\leq \frac{\varepsilon^{n}\|D\beta\|_{0}+\varepsilon^{n+m}\|Dr_{1}\|_{0}+C|\Delta\Omega|}
{1-\varepsilon^{n}\|D\beta\|_{0}\delta_{1}-\varepsilon^{n+m}\|Dr_{1}\|_{0}+C|\Delta\Omega|}\|w_2-w_1\|_0,
 \end{split}
\end{equation*}
\begin{equation*}
 \Bigg(\textrm{i.e.}\quad \textrm{Lip}_{2}v\leq  \frac{\varepsilon^{n}\|D\beta\|_{0}+\varepsilon^{n+m}\|Dr_{1}\|_{0}+C|\Delta\Omega|}
{1-\varepsilon^{n}\|D\beta\|_{0}\delta_{1}-\varepsilon^{n+m}\|Dr_{1}\|_{0}+C|\Delta\Omega|}
\Bigg)
\end{equation*}
and
\begin{equation*}
 \begin{split}
  |w_{2}(v(y_1,w_2))-w_1(v(y_1,w_1))|\leq &\|w_{2}-w_{1}\|_{0}+\|Dw_{1}\|_{0}\cdot|v(y_{1},w_2)-v(y_1,w_1)|\\
\leq& (1+\delta_{1} \textrm{Lip}_{2}v)\cdot \|w_2-w_1\|_{0}.
 \end{split}
\end{equation*}

Therefore, one easily gets
\begin{equation*}\label{contraction stable}
\begin{split}
 \|\mathscr{F}^{s}&[w^{s}_{2},w^{u}_{2}]-\mathscr{F}^{s}[w^{s}_{1},w^{u}_{1}]\|_{0}\\
\leq& \Bigg\{(1-C_{s}\varepsilon^{n})(1-\varepsilon^{n}\|D\beta\|_{0}\delta_{1})^{-1}+
\varepsilon^{n+m}\|Dr_{2}\|_{0}\cdot\Big[\textrm{Lip}_{2}v+(1+\delta_{1} \textrm{Lip}_{2}v)\Big]\\
&+2\varepsilon^{n}\delta_{0}\|\eta\|_{2}+\varepsilon^{n+1}\left(\|\eta\|_{1}+C\frac{|\Delta\Omega|}{\varepsilon^{n+1}}\right)\Bigg\}
\cdot \|(w_{2}^{s}-w_{1}^{s},w_{2}^{u}-w_{1}^{u})\|_{0},\\
\end{split}
\end{equation*}
which is also a contraction on $\mathscr{L}$ if $\varepsilon$ is chosen sufficiently small.

Together with Lemma \ref{invariant subset L} and the contractility of $\mathscr{F}$on $\mathscr{L}$, we obtain
from \cite[Lemma 2.4]{HdlL16} the following lemma on the existence and regularity of
the fixed point of $\mathscr{F}$.

\begin{lemma}\label{persistence of invariant surface}
Under the assumptions of Lemma \ref{invariant subset L}, there exists a  fixed point
$w^{*}\in \mathscr{C}^{L-1}(\TT^{r}, \RR^{d-r})$ of $\mathscr{F}$, which is unique in the closure of $\mathscr{L}$ in $\mathscr{C}^{0}(\TT^{r}, \RR^{d-r})$.
\end{lemma}

\noindent\textbf{Proof of Theorem \ref{phase locking theorem gen}.}
By the formulation of the functional $\mathscr{F}$, Theorem \ref{phase locking theorem gen}
is an immediate result of Theorem \ref{lemma normal form general} and Lemma \ref{persistence of invariant surface}.
The regularity $L-1$ of the invariant surface (given by the graph of $w^{*}$) is less than that of the resonant normal form.
\qed

\subsection{Lindstedt series for generic torus maps close to rotation}

Consider the following torus map
\begin{equation}
F_{\varepsilon}\begin{pmatrix}x\\y\end{pmatrix}=\begin{pmatrix}x+\omega+\varepsilon g(x,y)\\ y+\varepsilon h(x,y)\end{pmatrix} \in\TT^{r}\times\TT^{d-r},
\end{equation}
where $\omega\in\RR^{r}$ satisfies the Diophantine condition.

 We look for $l_{\varepsilon}=(l_{\varepsilon}^{x}, l_{\varepsilon}^{y}):\TT^{r}\rightarrow \TT^{d}$ and
$u_{\varepsilon}:\TT^{r}\rightarrow\TT^{r}$ such that
\begin{equation}\label{embedding general}
F_{\varepsilon}\circ l_{\varepsilon}=l_{\varepsilon}\circ u_{\varepsilon},
\end{equation}
by which the graph of $l_{\varepsilon}$ is the desired invariant surface.
To formulate the Lindstedt series, we
expand $l_{\varepsilon}$ and  $u_{\varepsilon}$ into power series in $\varepsilon$ as
\[
l_{\varepsilon}(\sigma)=\sum_{j\geq 0}\begin{pmatrix}l^{x}_{j}(\sigma)\\ l_{j}^{y}(\sigma)
 \end{pmatrix}
\varepsilon^{j},\quad
u_{\varepsilon}(\sigma)=\sum_{j\geq 0} u_{j}(\sigma)\varepsilon^{j},
\]
and then
take formal calculations.

By matching the coefficients of $\varepsilon^{0}$-terms, we obtain
$F_{0}\circ l_{0}=l_{0}\circ u_{0}$,
which can be solved by choosing
\[
 l_{0}^{x}(\sigma)=\sigma, \quad l_{0}^{y}(\sigma)=y_{0}=\textrm{Constant},\quad u_{0}(\sigma)=\sigma+\omega.
\]
For those $\varepsilon^{1}$-terms, we get
$$ F_{1}\circ l_{0}(\sigma)+DF_{0}\circ l_{0}(\sigma)\cdot l_{1}(\sigma)=l_{1}\circ u_{0}(\sigma)+Dl_{0}\circ u_{0}(\sigma)\cdot u_{1}(\sigma).$$
It follows that
\begin{equation}\label{Lindstedt 1st general}
 \left\{
\begin{aligned}
&l_{1}^{x}(\sigma+\omega)-l_{1}^{x}(\sigma)=g(\sigma,y_{0})-u_{1}(\sigma),\\
&l_{1}^{y}(\sigma+\omega)-l_{1}^{y}(\sigma)=h(\sigma,y_{0}).
\end{aligned}
\right.
\end{equation}
Assuming that there exists $y_{0}\in\TT^{d-r}$ such that
\begin{equation}\label{average vanishing general}
\langle h(\cdot, y_{0}) \rangle=\int_{\TT^{d-r}} h(\sigma, y_{0})\ \textrm{d}\sigma=0,
\end{equation}
then we can choose
\[
u_{1}(\sigma)\equiv \langle g(\cdot, y_{0})\rangle
\]
and solve \eqref{Lindstedt 1st general}
for $(l_{1}^{x},l_{1}^{y})$ by Lemma \ref{cohomology-torus} in KAM theory, leaving
the average
\begin{equation*}
  \langle l_{1}\rangle=(\langle l_{1}^{x}\rangle, \langle l_{1}^{x}\rangle)=\int_{\TT^{r}}
  l_{1}(\sigma)~\mathrm{d}\sigma
\end{equation*}
to be specified.

To clarify the induction, we proceed to compute the equation for those
$\varepsilon^{2}$-terms, which reads
$$ F_{2}\circ l_{0}+DF_{1}\circ l_{0}\cdot l_{1}+DF_{0}\circ l_{0}\cdot l_{2}=
 l_{2}\circ u_{0}+Dl_{1}\circ u_{0}\cdot u_{1}+Dl_{0}\circ u_{0}\cdot u_{2}.$$
It follows that
\begin{equation}\label{Lindstedt 2nd general}
 l_{2}(\sigma+\omega)-l_{2}(\sigma)=
\begin{pmatrix}
 Dg(\sigma, y_{0})l_{1}(\sigma)-Dl_{1}^{x}(\sigma+\omega)u_{1}-u_{2}(\sigma)\\
Dh(\sigma, y_{0})l_{1}(\sigma)-Dl_{1}^{y}(\sigma+\omega)u_{1}
\end{pmatrix}
.
\end{equation}
 Let $l_{1}(\sigma)=\langle l_{1}\rangle +\widetilde{l}_{1}(\sigma)$
and recall that $\widetilde{l}_{1}(\sigma)$ is uniquely determined by
 \eqref{Lindstedt 1st general}.
Then the R.H.S. of
\eqref{Lindstedt 2nd general} reads
\[
 \begin{pmatrix}
  Dg(\sigma, y_{0})\\ Dh(\sigma, y_{0})
 \end{pmatrix}
\cdot\widetilde{l}_{1}(\sigma)+
\begin{pmatrix}
Dg(\sigma, y_{0})\\ Dh(\sigma, y_{0})
\end{pmatrix}
\cdot\langle l_{1}\rangle-
\begin{pmatrix}
 Dl_{1}^{x}(\sigma)u_{1}+u_{2}(\sigma)\\ 0
\end{pmatrix}.
\]
In order the average of  R.H.S of \eqref{Lindstedt 2nd general} to vanish, we need to choose
parameters $\langle l_{1}\rangle=(\langle l_{1}^{x}\rangle, \langle l_{2}^{y}\rangle)$
and $u_{2}(\sigma)$. More precisely, we have
\begin{equation*}
 \begin{split}
&\langle D_{2}h(\cdot, y_{0})\rangle\cdot \langle l_{1}^{y}\rangle= \langle
Dh(\cdot, y_{0})\widetilde{l}_{1} \rangle,\\
 &\langle D_{2}g(\cdot, y_{0})\rangle\cdot \langle l_{1}^{y}\rangle
-\langle u_{2}\rangle=
\langle Dg(\cdot, y_{0})\widetilde{l}_{1} \rangle,
 \end{split}
\end{equation*}
since $\langle D_{1}g(\cdot, y_{0})\rangle=
\langle D_{1}h(\cdot, y_{0})\rangle=\langle Dl_{1}^{x}\rangle=0
$.

Then if
\begin{equation}\label{average non-degeneracy general}
 \det \langle D_{2}h(\cdot, y_{0})\rangle\neq 0,
\end{equation}
we just take
\begin{equation*}
  \langle l_{1}^{y}\rangle=\langle  D_{2}h(\cdot, y_{0})\rangle^{-1}\cdot
\langle
Dh(\cdot, y_{0})\widetilde{l}_{1} \rangle ,
\end{equation*}
and
\begin{equation*}
u_{2}(\sigma)\equiv \langle  D_{2}g(\cdot, y_{0})\rangle\cdot\langle  D_{2}h(\cdot, y_{0})\rangle^{-1}\cdot
\langle
Dh(\cdot, y_{0})\widetilde{l}_{1} \rangle-
\langle Dg(\cdot, y_{0})\widetilde{l}_{1} \rangle.
\end{equation*}

Now the induction is clear. By similar arguments in Section 2.2, we obtain
the following result.

\begin{proposition}\label{appen Lind}
 Let $\omega\in\RR^{r}$ satisfy Diophantine condition. Assume there exists $y_{0}\in\RR^{d-r}$ such
that \eqref{average vanishing general} and
\eqref{average non-degeneracy general} hold. Then we can find two formal power series
\[
 l_{\varepsilon}(\sigma)=\sum_{j\geq 0} l_{j}(\sigma)\varepsilon^{j}\quad\textrm{and}\quad
u_{\varepsilon}(\sigma)=\sigma+\omega+\sum_{j\geq 1} u_{j}\varepsilon^{j}
\]
such that \eqref{embedding general} holds, where $\{u_{j}\}_{j\geq 1}$ are
a sequence of constant vectors.

\end{proposition}

We are also able to establish the relationship
between the Lindstedt series and the true embedding of the invariant surface as that in
Theorem \ref{Lind theorem}.


\def\cprime{$'$} \def\cprime{$'$}

\end{document}